\documentclass[12pt]{article}

\usepackage{amsmath,amsthm} 
\usepackage{amssymb}    
\usepackage{graphicx}     
\usepackage{hyperref} 
\usepackage{url}
\usepackage{amsfonts} 

\newcommand{\comp}{\;\raisebox{2pt}{$\scriptstyle \circ$} \;}

\newcommand{\bbRR}{\mbox{${\mathbb R}^2 $}}
\newcommand{\bbC}{\mbox{$\mathbb C $}}

\newcommand{\figref}[1]{{\textsc Figure}~\ref{#1}}

\newcommand{\colvv}[2]
     {\left[ 
         \begin{array}{c} #1 \\ #2  \end{array}
      \right]}

\newcommand{\matrrvv}[4]
     {\left[ 
         \begin{array}{cc} #1 & #2 \\ 
                           #3 & #4 \end{array} \right]}

\newtheorem{theorem}{Theorem}[section]

\newtheorem{lemma}{Lemma}[section]

\newtheorem{definition}{Definition}[section]

\begin{document}

\title{Currie's Mysterious Pattern and Iterated Functions}

\author{Dan Kalman \\
American University \\
Washington, DC 20016 \\
\texttt{kalman@american.edu}}

\maketitle

\begin{abstract}
In \cite{mathr&r}, Currie discusses what he calls a {\em
mysterious pattern} involving the sequence 
\[ a_{n} = 2^n \sqrt{2 - \sqrt{2 + \sqrt{2 + \cdots + \sqrt{2}}}},\] 
where $n$ is the number of nested radicals. 
Part of the mystery is that $a_n$ converges to $\pi.$

In this paper we discuss a general framework for results like the
mysterious pattern in the context of iterated functions.
\end{abstract}

In \cite{mathr&r}, Currie discusses what he calls a {\em
mysterious pattern} involving the sequence
\begin{equation}  \label{mp1}
a_{n} = 2^n \sqrt{2 - \sqrt{2 + \sqrt{2 + \cdots + \sqrt{2}}}},
\end{equation}
where $n$ is the number of nested radicals.  Part of the mystery is
that $a_n$ converges to $\pi.$

It is natural to ask whether this is an isolated result, or whether it
is part of some larger family of related patterns and limits.  And in
fact, it is somewhat surprising how many paths start with the
mysterious pattern and lead to interesting
extensions, insights, or generalizations.  For example, nested radical
expressions similar to those in~\eqref{mp1} but with different sign
patterns have been considered in several works, starting
with~\cite{polyaszego}.  In~\cite{geom_ext_ppr}
several extensions are discussed based on the original geometric context of Currie's 
pattern.  Taking a different
direction, $a_n$ can be viewed as one instance of an entire family of
algebraically similar sequences.  In particular, fixing a real
parameter $x$, the sequence  
\[ a_{n}(x)=  
  (2x)^{n/2} \sqrt{x - \sqrt{x^2-x + \sqrt{x^2-x + \cdots + \sqrt{x^2-x}}}} \]
converges for all $x > 1$, and reduces to $a_n$ when $x=2$.  This
leads naturally to the function $C(x) = \lim_{n\rightarrow \infty}
a_n(x)$, the study of which reveals interesting results and
conjectures as detailed in~\cite{curriefunc}.

As Currie shows in~\cite{mathr&r}, his mysterious pattern admits a formulation as an
infinite product, equivalent to Vi\`{e}te's 1593 product formula for
$2/\pi$, cited in~\cite{levin-ramijrnl}.  There are quite a few papers
that have considered variations of Vi\`{e}te's identity, for example 
\cite{garcia-moreno-prophet-FQ}, \cite{levin-ramijrnl},
\cite{moreno-garcia-mathmag}, \cite{moreno-garcia-jrnl-approx-thry},
\cite{osler-fibquart}, and \cite{osler-mathgaz}.  In particular,
\cite{moreno-garcia-jrnl-approx-thry} relates one of its results to the
nested radical sequences of~\cite{polyaszego}.  But we have found no
other instances where a generalized product identity has been used to
derive an equivalent result in a nested radical form analogous to~\eqref{mp1}.

In the present paper we relate Currie's mysterious pattern to a
broader class of sequences and limits 
defined in terms of iterated functions.  These take the form
\begin{equation}  \label{itfunc1}
a_{n} = (1/r^n)|L-f^{(n)}(t_0)|. 
\end{equation}
Here $f^{(n)}(x)$ represents the $n$th iterate of a function $f,$ so
$f^{(2)}(x)=f(f(x))$, $f^{(3)}(x) = f(f(f(x)))$, and so
on. For the cases we will be interested in, 
$L = \lim_{n \rightarrow \infty} f^{(n)}(t_0)$, and $1/r$ is a positive
factor for which $a_n$ converges to a finite positive limit.
Observe that squaring \eqref{mp1} yields an instance of \eqref{itfunc1}
when $f(t)= \sqrt{2+t}$, $L=2$, and $r=1/4$.  We will consider conditions on
$f$ for which suitable $L$ and $r$ exist, and the possibility of
evaluating the limit of the corresponding sequence $a_n$.  

\section{Iterated Function Basics}
We consider real functions $f$ defined on a closed interval $I$, and smooth
enough to permit us to differentiate as many times as necessary.  For
a given initial value $t_0 \in I$, the sequence
$\{f^{(n)}(t_0)\}$ is sometimes
referred to as the {\em orbit} of $t_0$ under $f,$ calling to mind a
dynamical perspective.  

The dynamics of function iteration have been widely studied.  In
general terms, if $f$ is a contraction map on a complete metric space
$D$, then there is a unique fixed point $L$ in $D$ and from any
initial $t_0$ in $D$ the sequence $\{f^{(n)}(t_0)\}$ converges to
$L$.  We are interested in the special case of a smooth $f$ mapping
$I$ into itself.  If for some fixed $k$,
$|f'(t)| \le k < 1$ on $I$, then for any $t$ and $s$ in $I$, $|f(t) -
f(s)| \le k|t-s|$ so $f$ is a contraction with fixed point $L$.  We
assume that $L$ is in the interior of $I$, which we can verify in the
contexts of interest for this paper, and denote $|f'(L)| = m \le k.$
Then for arbitrarily small positive $\epsilon$, $|f'(t)|$ must fall
between $m-\epsilon$ and $m + \epsilon$ for all $t$ in some
neighborhood $N$ of $L$.  Now if $t \in N$ the mean value theorem implies
$|L - f(t)| = |f(L) - f(t)| = |f'(\xi)(L - t)|$ for some $\xi$
between $L$ and $t$, and hence in $N$.  Therefore
$(m-\epsilon)|L-t| \le |L - f(t)| \le (m+\epsilon)|L-t|$.  
Applying this observation repeatedly, we find
\begin{equation} \label{inequ1}
(m-\epsilon)^n|L-t| \le |L - f^{(n)}(t)| \le (m+\epsilon)^n|L-t|
\end{equation}
for all $n \ge 1$ and all $t$ in $N.$   This shows that iterating $f$
from any starting point in $I$ will eventually converge to $L$
essentially exponentially.

To relate the mysterious pattern to function iteration we
make two minor algebraic modifications in~\eqref{mp1}.  First we square both
sides to define
\[ c_{n} = 
  4^n \left(2 - \sqrt{2 + \sqrt{2 + \cdots + \sqrt{2}}}\right),\]
which we know has limit $\pi^2$ (and only $n-1$ radicals).  Then we
rewrite $c_n$ in the form 
\begin{equation}  \label{mp2}
c_{n} = \frac{2 - \sqrt{2 + \sqrt{2 + \cdots + \sqrt{2}}}}
                 {(1/4)^n}.
\end{equation}
In the numerator we see an iterated function subtracted from 2; in the
denominator we have an exponential sequence that decays to zero.  And
since we know that $c_n$ has a finite positive limit, we infer that the
numerator is also decaying to zero asymptotically as $(1/4)^n$.  Thus,
defining $f(t)=\sqrt{2+t}$, we see that $f^{(n)}(0)$ has limit 2 and
\begin{equation}  \label{mp3}
\lim_{n \rightarrow \infty} \frac{2 - f^{(n-1)}(0)}{(1/4)^n}=\pi^2.
\end{equation}
An equivalent form that is sometimes more convenient is 
\begin{equation}  \label{mp3.1}
\lim_{n \rightarrow \infty} \frac{2 - f^{(n)}(0)}{(1/4)^n}=\frac{\pi^2}{4}.
\end{equation}

As an obvious generalization, we can define a sequence
\begin{equation} \label{obvgen}
 c_n = \frac{L - f^{(n)}(t_0)}{r^n}
\end{equation}
for some function $f$ and with positive constants $L$ and $r$.
When $r < 1,$  $c_n$ can only converge when $f^{(n)}(t_0)$ converges to
$L$, which is thus a fixed point of $f$.  

What can we say about $r$?  From~\eqref{inequ1} we know that  
\begin{equation} \label{convbnds1}
   \left(\frac{m-\epsilon}{r}\right)^n|L-t_0| 
        \le \frac{|L - f^{(n)}(t_0)|}{r^n} 
        \le \left(\frac{m+\epsilon}{r}\right)^n|L-t_0|.
\end{equation}
This shows that $c_n$ has a finite positive limit only if $r=m.$  For
if $r < m$ we can choose 
$\epsilon$ so that $r < m-\epsilon$, driving $|L-f^{(n)}(t_0)|/r^n$ to
infinity.  Likewise, if $r > m$, choosing $\epsilon$ so that $r >
m+\epsilon$ drives $|L-f^{(n)}(t_0)|/r^n$ to 0.

Although $r = m$ is necessary for $c_n$ to have a finite positive
limit, it is not obvious whether it is sufficient.  With $r = m$,
\eqref{convbnds1} becomes 
\begin{equation} \nonumber
   \left(1 - \frac{\epsilon}{m}\right)^n|L-t| 
        \le \frac{|L - f^{(n)}(t)|}{m^n} 
        \le \left(1+\frac{\epsilon}{m}\right)^n|L-t|,
\end{equation}
and in the limit we obtain no information about the convergence of the
ratio in the middle.  On the other hand, we have not discovered an
example demonstrating that $r=m$ is not sufficient for 
convergence of $c_n$ to a finite positive limit.

In the case of the mysterious pattern as expressed in \eqref{mp2},
$f(t) = \sqrt{2+t}$ for which $0 < f'(t) < 1/2$ on $[-1,\infty)$.  At
the unique fixed point $L = 2$, we have $f'(L) = 1/4$, making $r =
1/4$.  This shows that the convergence of the sequence in~\eqref{mp2}
to $\pi^2$ represents convergence of an instance of the more
general case defined in~\eqref{obvgen}.

\subsection{Candidate Sequences and Proper Convergence}

This leads to a general class of sequences generalizing the one in
\eqref{mp2}. For a given
contraction $f$ with fixed point $L$, we consider sequences of the
form  
\begin{equation} \label{resid seq 1}
  c_n =  \frac{|L - f^{(n)}(t_0)|}{|f'(L)|^n},
\end{equation}
which we will call {\em candidate sequences}.  A candidate sequence
with a nonzero finite limit will be said to converge {\em properly}.  
Note that taking $t_0 = L$ leads to $c_n
= 0$  for all $n$.  In this case $c_n$ does converge, but not
properly.

Borrowing some terminology from the area of differential equations, we
can think of $f^{(n)}(t_0)$ as samples of a signal that is approaching
a steady state of $L$, and consider $|L - f^{(n)}(t_0)|$ as a
transient component.  When $c_n$ converges properly, the transient
decays approximately exponentially, and $c_n$ is the residual when we
divide out the exponential component of that decay.  This residual
depends on the value of $t_0$ as it accumulates throughout the decay process.
The closer $t_0$ is to $L$, the less we expect the residual to be.

This intuition is confirmed by considering what happens when we
introduce a delay in the sequence of iterates.  For example, consider
a candidate sequence 
\[ c_n^{(5)} = \frac{|L - f^{(n)}(t_5)|}{|f'(L)|^n}\]
where $t_5 = f^{(5)}(t_0).$  Then we have
\[ c_n^{(5)} = \frac{|L - f^{(n+5)}(t_0)|}{|f'(L)|^n} 
             =|f'(L)|^5 \cdot \frac{|L - f^{(n+5)}(t_0)|}{|f'(L)|^{n+5}} 
             =|f'(L)|^5 c_{n+5} .\]

Thus, if $c_n$ has limit $M$, then $c_n^{(5)}$ has limit $|f'(L)|^5M.$
Here $c_n$ and $c_n^{(5)}$ are both candidate sequences for the
same function $f$, and differ only in the initial value from which $f$
is iterated, yet they have different residuals.  Moreover, $t_5$ is 
closer to $L$ than $t_0$, and the residual for $c_n^{(5)}$ is less
than the residual for $c_n$.

The class of properly converging
candidate sequences generalizes the case of the mysterious pattern.
As an example, we let $f(t) = \sqrt[3]{24+t},$ which has a fixed
point at $L = 3.$  We find $|f'(L)| = 1/27$.  The corresponding
candidate sequence is
\[   c_n =  27^n |3 - f^{(n)}(t_0)| . \]
To strengthen the analogy with the mysterious pattern, it makes sense to
look at the cuberoot of $c_n$, thus defining
\[ a_n = 3^n \sqrt[3]{3-f^{(n)}(t_0)}.   \]
Then
\[ a_n = 3^n\sqrt[3]{3-\sqrt[3]{24 + \sqrt[3]{24 + \cdots \sqrt[3]{24+t_0}}}}  \]
where the number of radicals is $n+1$.  Later we will show that this
sequence does converge properly, providing a nice analog of the
mysterious pattern.  However, we have not discovered the exact value
of the limit.

Generalizing the prior example, for a positive $L$ we define
$f(t) = \sqrt[3]{L^3 - L + t},$ which evidently has a fixed point at
$t = L.$  Then for $t > L - L^3$
\[ f'(t) = \frac{1}{3(L^3 - L + t)^{2/3}}    \]
so that $f'(L) = 1/(3L^2).$  For any $L > \sqrt{1/3}$ we have
$0 < f'(L) < 1$.  In this case, $f$ is a contraction in the interval 
$(L - L^3 + 1/\sqrt{27},\infty)$, with $L$ the unique fixed point in
that interval.  Taking any $t_0$ in the interval, we obtain a
candidate sequence
\[ c_n = \frac{\left|L - f^{(n)}(t_0)\right|}{f'(L)^n} 
       = (3L^2)^n\left|L - f^{(n)}(t_0)\right|.   \]
As in the preceding example, $\sqrt[3]c_n$ is an analog of the
mysterious pattern sequence for any $L > \sqrt{1/3}$.

Returning to the more general context, we should like to characterize
the candidate sequences that converge properly.  We have obtained a
partial result in this direction, identifying a class of functions
that have properly convergent candidate sequences.
This class of functions includes our
original example $f(t) = \sqrt{2+t}$, as well as the preceding
parameterized example $f(t) = \sqrt[3]{L^3 - L + t}$.  For such
functions we demonstrate proper convergence of candidate sequences by
showing that they are monotonic and bounded above and away from 0.  To obtain
the necessary bounds, we compare a candidate sequence for $f$ with
corresponding sequences for functions that closely approximate $f$
near $L.$  

Next we will derive the monotonicity result.  Later we will discuss an
approximation to $f$ that can be used to obtain a necessary bound.

\subsection{Monotonicity of Candidate Sequences}
Throughout this discussion, we assume $f$ is a contraction on an
interval $I,$ with fixed point $L$ in the interior of $I$, and that $f'(t) > 0$ and
$f''(t)<0$ for $t \in I$.  We denote $f'(L)$ as $m$ and observe $0 < m
< 1$.  As a first step, we show that sequences of the form
$f^{(n)}(t_0)$ are monotonic for all $t_0 \ne L$ in $I.$

\begin{lemma}  \label{lemma1}
Let $t_0 \ne L$ be an element of $I$.  Define $t_n =
f^{(n)}(t_0)$. Then $t_0 > L$ implies $t_n$ is decreasing with
limit $L$ and $t_0 < L$ implies $t_n$ is increasing with limit $L$. 
\end{lemma}

\begin{proof}
For $t \in I$ and $t > L$ we have
$f(t) - L = f(t) - f(L) = f'(\xi)(t - L)$ for some $\xi$ in $I$.
From the assumed properties of $f$, $0 < f'(\xi) < 1$.  So $0 < f(t) -
L < t - L$ and $L < f(t) < t$.  For any $t_0 > 
L$ these results imply that
$t_n$ is strictly decreasing and bounded below by $L$.  Consequently
$t_n$ must converge, and the limit, being a fixed point of $f,$ must
equal $L$.  This proves the first part of the lemma.  The proof of the
second part is similar. 
\end{proof}


Next we prove that candidate sequences are also monotonic.

\begin{lemma}   \label{lemma2}
Let $t_0 \ne L $ be in $I$ and with 
$t_n = f^{(n)}(t_0)$ define the
candidate sequence  $c_n = |L-t_n|/m^n$ .  
Then $t_0 > L$ implies $c_n$ decreases with $n$ and
$t_0 < L$  implies $c_n$ increases with $n$.
\end{lemma}

\begin{proof}  Suppose  
$t_0 > L$, then each $t_n > L$ so $|L - t_n| = t_n - L$.  Now for $n
\ge 0,$ we have $t_{n+1} - L = f(t_n) - f(L) = f'(\xi)(t_n-L)$ for
some $\xi$ between $L$ and $t_n$.  Thus $f'(\xi) \le f'(L) = m$,
which implies that $t_{n+1} - L \le m(t_n-L)$.  Now divide by
$m^{n+1}$ to conclude $(t_{n+1} - L)/m^{n+1} \le (t_n-L)/m^n$.  This
shows that the candidate sequence $c_n$ decreases
with $n$ for $t_0 > L.$  

The proof of the $t_0 < L$ case is similar.

\end{proof}

The lemma establishes the monotonicity of candidate sequences.  For
the case of an initial term less than $L$, with an increasing
candidate sequence, we obtain proper
convergence by exhibiting a finite upper bound.
Similarly, when the initial term is greater than $L$, we obtain proper
convergence by exhibiting a positive lower bound.  As mentioned, we
obtain such bounds by comparing a contraction $f$ with a closely
related approximating function.  We
proceed shortly to the discussion of such functions.  But first, we
derive two additional useful results.

\begin{lemma}  \label{lemma3}
With the restriction that $t_0 \in I,$ if the candidate sequence for some $t_0 < L$
converges properly, then the candidate sequence for every $t_0 < L$
converges properly, and if the candidate sequence for some $t_0 > L$
converges properly, then the candidate sequence for every $t_0 > L$
converges properly.  
\end{lemma} 

\begin{proof}  Suppose the candidate sequence for $t_0 < L$
converges properly.  Consider $s_0$ such that 
$t_0 < s_0 < L$.  We denote $f^{(n)}(t_0)$ as $t_n$ and $f^{(n)}(s_0)$
as $s_n$.  We know both $s_n$ and $t_n$ increase to a limit of $L$.
To show that the candidate sequence $(L-s_n)/m^n$ converges
properly, it suffices to show it is bounded above.

Since $f$ is increasing, $t_0 < s_0$ implies $t_1 = f(t_0) < f(s_0) =
s_1$ and by induction, $t_n < s_n$ for all $n$.  Therefore
$L-s_n < L-t_n$, hence
\[ \frac{L-s_n}{m^n} < \frac{L-t_n}{m^n}.       \]
We have assumed that the sequence on the right converges properly, and
hence it is bounded above.  This shows that the candidate sequence
$(L-s_n)/m^n$ is bounded above, as desired.

On the other hand, if $s_0 < t_0,$ since $s_n$ increases to $L$ there
exists some $k$ such that  
$t_0 < s_k$.  From the prior case we know that the candidate sequence
$(L - f^{(n)}(s_k))/m^n$ converges properly.  But $f^{(n)}(s_k) =
s_{k+n}$.  So we know that
\[ \frac{L-s_{k+n}}{m^n} = m^k \frac{L-s_{k+n}}{m^{k+n}}    \]
converges to a finite positive limit.  Therefore so does
$(L-s_n)/m^n.$  This completes the proof of the first part of the
lemma.

For the second part suppose $t_0 > L$ and the corresponding candidate
sequence converges properly.  Again we consider $s_0 \ne t_0$ and
define $t_k$ and $s_k$ as before.  However for this part of the proof
the candidate sequences associated with $t_n$ and $s_n$ are
decreasing, and it suffices to show that $s_n$ has a positive lower
bound.

If $s_0 > t_0,$ then $s_n > t_n$ for all $n$ because $f$ is
increasing.  Thus $s_n - L > t_n - L$ and hence
\[ \frac{s_n - L}{m^n} > \frac{t_n - L}{m^n}.   \]
By assumption, the sequence on the right has a positive limit $M$ which is
also a lower bound.  Thus $M$ is also a positive lower bound for the
sequence on the left.  This shows that the candidate sequence for
$\{s_n\}$ is properly convergent.

On the other hand, suppose $L < s_0 < t_0$.  We know that $t_n$
decreases to $L$ so for some $k,$ $L < t_k < s_0.$  Moreover, the
sequence $f^{(n)}(t_k)$ has a properly convergent candidate sequence
because
\[ \frac{f^{(n)}(t_k) - L}{m^n} =
   \frac{t_{k+n}-L}{m^n} =
   m^k \frac{t_{k+n}-L}{m^{k+n}}.  \]

\noindent
This shows the candidate sequence for $f^{(n)}(t_k)$
is a constant positive multiple of a subsequence of the candidate
sequence for $t_n$, which we know converges properly.  And now we have
reduced this case to that of the preceding paragraph. We have the
initial point $t_k$ for a sequence of iterates under $f$ having a
properly convergent candidate sequence, and a greater initial point
$s_0$.  Therefore, the candidate sequence for $\{s_n\}$ must also
converge properly. 
\end{proof}

The final lemma is useful in comparing candidate sequences for two
different functions.

\begin{lemma}   \label{lemma4}
Suppose $f$ and $g$ are both contractions on an
interval $I$ with fixed point $L \in I$, that $f'(L)=g'(L)=m<1,$ and
that both $f$ and $g$ are increasing in
$I$.  For $t_0$ and $u_0$
in $I$ define sequences
$\{t_n\} = \{f^{(n)}(t_0)\}$ and $\{u_n\} =
\{g^{(n)}(u_0)\}$. We also have candidate sequences defined by 
$c_{f:n} = |L-t_n|/m^n$ and $c_{g:n} = |L-u_n|/m^n$.

Then the following hold for $n \ge 1$:

\hspace*{.3in}
\begin{minipage}{.8\textwidth}
\begin{itemize}
\item[(i).] if $g(t) < f(t)$ for $t<L$, and if $u_0 \le t_0 < L, $ then
$c_{f:n} < c_{g:n}.$
\item[(ii).] if $g(t) < f(t)$ for $t>L$, and if $t_0 \ge u_0 > L, $ then
$c_{f:n} > c_{g:n}$.
\end{itemize}
\end{minipage}

\end{lemma}

\begin{proof}
Suppose first that $g(t) < f(t)$ for $t < L$ and that $u_0 \le t_0 <
L.$  Then $u_1 = g(u_0) \le g(t_0)$ because $g$ is increasing, and 
$g(t_0) < f(t_0) = t_1.$ Proceeding by induction, if
$u_n < t_n$ then $u_{n+1} = g(u_n) < g(t_n) < f(t_n) = t_{n+1}$.  Thus
$u_{n+1} < t_{n+1}$.  This shows 
that $u_n < t_n$ for all $n \ge 1.$  Continuing, $L - t_n < L - u_n,$
and hence $(L - t_n)/m^n < (L - u_n)/m^n.$  But this shows that 
$c_{f:n} < c_{g:n}$ for $n \ge 1,$ because $t_n \le L$ and $u_n
\le L$ for all $n.$  This establishes (i).

Next suppose that $g(t) < f(t)$ for $t > L$ and that $t_0 \ge u_0 >
L.$  Then, arguing similarly as before, we can show that $t_n > u_n$
for all $n \ge 1.$  Proceeding, $t_n - L > u_n - L,$ 
and hence $(t_n - L)/m^n > (u_n - L)/m^n.$  But this shows that 
$c_{f:n} > c_{g:n}$ for $n \ge 1,$ because $t_n \ge L$ and $u_n
\ge L$ for all $n.$  This establishes (ii).
\end{proof}

\section{Root-Like Functions}
We are now nearly in a position to demonstrate convergence of
candidate sequences for a broad class of functions, characterized in
the following definition.

\begin{definition} \label{rootlikedef}
A smooth function $f$ is {\em root-like} in an interval $I$ when
\newline \rule{.2in}{0in}\parbox{.75\textwidth}
{\begin{enumerate}
\item[i.] $f$ is a contraction on $I,$
\item[ii.] $f$ has a fixed point $L$ in the interior of $I,$ 
\item[iii.] $f'(t) > 0$ and $f''(t) < 0$ in $I$.
\end{enumerate}}
\newline
We denote $f'(L)$ as $m$ and observe $0 < m < 1$.  
\end{definition}
Such a function is root-like in the sense that its graph near $L$ is similar
to a curve of the form $y = \sqrt[m]{x+H} - K$ near its (greatest) fixed point.
Requiring $L$ to be an interior point of $I$ excludes cases where the
fixed point is a double root of $f(x) = x$, as happens for example
with $f(x) = \sqrt{x-1/4}$.  The results about root-like functions to
follow can be extended to this excluded case by similar analyses, but
in the interest of brevity we leave the details to the interested reader.
In a similar spirit, for the sake of simplicity we assume $L>0.$

For a given $t_0 \in I,$ the
sequence $\{t_n\}$ is defined by $t_n = f^{(n)}(t_0)$.  We have seen
that such a sequence is monotonic, as is the corresponding candidate
sequence.  And we have seen that convergence of a candidate sequence
for one function can imply convergence of a candidate sequence for a
related function.

The next step is to introduce a function $Q$ that approximates $f$
near $L$ and for which the candidate sequence converges properly.  To
that end, we state the following definitions.

\begin{definition} \label{Q func def0}
A $Q$-function is a root-like function of the form
\begin{equation} \label{Qdef1}
Q(t) = \frac{at+b}{t+d} 
\end{equation}
where $a, b,$ and $d$, are constants.
\end{definition}

\begin{definition} \label{Q func def1}
Let $f$ be a root-like function with fixed point $L$ where $f'(L)=m
\in (0,1)$. We say $Q$ is an
{\em associated Q-function} for $f$ when $Q$ is a $Q$-function and
\newline \rule{.2in}{0in}\parbox{.5\textwidth}
{\begin{enumerate}
\item[i.] $Q(L) = f(L) = L$,
\item[ii.] $Q'(L) = f'(L) = m$, and
\item[iii.] $Q''(L) < f''(L)$.
\end{enumerate}}
\end{definition}

Complex functions of the form in~\eqref{Qdef1} are referred to as M\"obius
transformations, about which there is an extensive literature
with applications to several areas of mathematics.  Although we will
use some of the properties of M\"obius transformations, we will do so only
in the real variable context.  Hence we will refer to $Q(t)$
as a M\"obius {\em function}.  Note that this conflicts with the more
usual meaning of a family of functions $\mu(n)$ in number theory and combinatorics.

The M\"obius function $Q(t)$ is associated with the matrix
\[ M_Q = \matrrvv{a}{b}{1}{d},   \]
which leads to important insights about $Q$.  This will be discussed in
detail in Section~\ref{mobius sxn}.  For now we will restrict our
attention to a few properties 
that will be useful in showing convergence of candidate sequences.

As a preliminary step, let us see how $Q-$functions can be
constructed.  Assume that $Q(t)$ has a fixed point at $L$, that $Q'(L)
= m$, and that $Q''(L) = s$.  Direct computation leads to 
\begin{eqnarray}
\frac{aL+b}{L+d}       & = & L   \nonumber \\
\frac{ad-b}{(L+d)^2}   & = & m   \label{abd to mLs}\\
\frac{2(b-ad)}{(L+d)^3}   & = & s.  \nonumber
\end{eqnarray}
This system can be inverted to express $a, b,$ and $d$ as functions of
$L, m,$ and $s$.  If  $m \ne 0$ and $s
\ne 0$ then the values of $a,
b,$ and $d$ are given by
\begin{eqnarray}
a   & = & \frac{Ls-2m^2}{s}     \nonumber  \\
b   & = & \frac{L}{s}(2m^2-2m-Ls) \label{mLs to abd}\\
d   & = & -\frac{2m+Ls}{s}. \nonumber
\end{eqnarray}
These equations show that we can define an entire family of
$Q$-functions associated with a specified root-like $f.$  We have to
choose $L$ and $m$ as dictated by conditions i and ii in
Definition~\ref{Q func def1}, but we can choose any real $s < f''(L)$.  And for whatever
value of $s$ we select, the result will be a $Q$-function associated
with $f.$  The next lemma shows that all of these $Q$ functions share
an important characteristic.

\begin{lemma}  \label{distinct-eigs-lemma}
If $Q$ is an associated $Q$-function for a root-like
function $f,$ then the matrix $M_Q$ has distinct real eigenvalues.
\end{lemma}

\begin{proof}
The eigenvalues of $M_Q$ are real and distinct if and only if the
characteristic polynomial of $M_Q$ has positive discriminant, and
hence when $(a-d)^2+4b > 0.$  In terms of the parameters $L, m,$ and
$s = Q''(L)$, this condition becomes
\[ \left( \frac{2Ls+2(m-m^2)}{s} \right)^2 
+ 4\frac{L}{s}(2m^2-2m-Ls) > 0.     \]
We introduce $\mu = m-m^2$ and further simplify to
\[ 4\frac{(Ls+\mu)^2}{s^2} 
- 4\frac{Ls(2\mu+Ls)}{s^2}> 0.     \]
Dividing out the positive factor $4/s^2$ we reach the equivalent
inequality
\[ (Ls+\mu)^2 - 2\mu Ls - L^2s^2 > 0.   \]
This is evidently true because $0<m<1$, thus completing the proof.
\end{proof}

The significance of distinct eigenvalues is shown by the next lemma.

\begin{lemma} \label{c_n conv for Q-func}
Candidate sequences for $Q$-functions converge properly.
\end{lemma}

This conclusion follows from aspects of M\"obius functions that will
be discussed in the next section.  The key idea is that $M_Q$ is
diagonalizable because it has distinct eigenvalues.  That leads to an
explicit formula for $Q^{(n)}(t_0)$ as a function of $n$.
Consequently the limit of the corresponding candidate sequence can be
directly evaluated.  In particular, the limits are positive and
finite.

We can say more about the convergence of $Q$'s candidate
sequence $c_n$. Notice that $Q$ is itself root-like in some interval
$I$. Applying Lemma~\ref{lemma2} to $Q$ shows that for $t_0 \in I,$
$c_n$ increases if $t_0 < L$, and it decreases if $t_0 > L$.  In
the former case, $c_n$ is bounded above by its limit which is finite;
in the latter case it is bounded below by its limit which is
positive. 

This brings us to the following theorem. 

\begin{theorem} \label{thm1}
Let $f$ be root-like in $I$, with fixed point $L$ and $f'(L) = m$.
Then, for any $t_0 \in I,$ $t_0 \ne L$, 
the candidate sequence $|L-f^{(n)}(t_0)|/m^n$ converges properly.
\end{theorem}

\begin{proof}
Note that $0 < m < 1.$  Let $Q$ be a $Q$-function associated with $f$.
Then $Q$ is root-like on some interval $I_1$ containing $L$.  We know that
$f''(L) > Q''(L)$, so $(f-Q)''(t) > 0$ in some interval $I_2$ containing
$L$.  Define $N = I \cap I_1 \cap I_2$. 

By assumption, $f-Q$ and its first derivative vanish at $L$, and
$(f-Q)''$ is positive in $N$.  By Taylor's theorem, we have
\[ (f-Q)(t) = \frac12 (f-Q)''(\xi)(t-L)^2  \]
for some $\xi$ between $t$ and $L$.  The right side of this equation
is evidently positive, so $Q(t) < f(t)$ for $t \in N.$

Now we can apply Lemma~\ref{lemma4} comparing candidate sequences for $f$ 
and $Q$.  Following the notation of the lemma, define sequences $t_n =
f^{(n)}(t_0)$ and $u_n = Q^{(n)}(u_0)$ with initial terms $t_0=u_0 \ne
L$ in $N$. Denote the corresponding candidate sequences by $c_{f:n}$ and
$c_{Q:n}$, respectively.  

There are two cases.  If $t_0 = u_0 < L$, then 
$c_{f:n} < c_{Q:n}$.  By Lemma~\ref{lemma2}, $c_{f:n}$ is increasing.  By earlier
remarks $c_{Q:n}$ is bounded above by its limit.  Thus
$c_{f:n}$ is increasing and bounded above and so converges to a finite
limit.  Furthermore we know the limit is positive because $c_{f:n} \ge
0$ by definition and also is increasing to its limit.

Similarly, if $t_0 = u_0 > L$, Lemma~\ref{lemma4} shows that 
$c_{f:n} > c_{Q:n}$.  Now $c_{Q:n}$ is decreasing and is bounded below
by its limit, which is positive.  Thus
$c_{f:n}$ is decreasing and bounded below and so converges to a finite
limit.  And the limit is positive because the lower bound is positive.

So far we have shown that $c_{f:n}$ converges properly for any $t_0
\in N,$ $t_0 \ne L$.  By Lemma~\ref{lemma3}, the same conclusion must hold with
$I$ in place of $N.$  This completes the proof.
\end{proof}

\subsection*{Examples}  Earlier we considered
the functions $f(t) = \sqrt[3]{L^3 - L + t}$ as a
generalization of the function $f(t) = \sqrt{2+t}$ that arises in
connection with the original mysterious pattern.  We can generalize
even further by defining $f(t) = \sqrt[k]{L^k - L + t}$ for $L>1$ and
integer $k \ge 2.$  Notice that $f(L) = L$.  Also, 
$f'(t) = \frac{1}{k}\left(L^{k}-L+t\right)^{\frac{1}{k}-1}$ so
$m = f'(L) = \frac{1}{kL^{k-1}} \le \frac{1}{2L} \le \frac{1}{2}.$  In
fact, $f'(t)$ is less than $1$ for 
$t > L-L^{k}+\left(\frac{1}{k}\right)^{\frac{k}{k-1}}$.  This defines
the maximal interval in which $f$ is a contraction.  As a
simplification, if we assume that $L > (1+\sqrt{2})/2 \approx 1.21,$
$f$ will be a contraction on $(0,\infty)$.


For this family of functions, candidate sequences always converge
properly.  By way of verification, it is sufficient to check that 
$f''(L) < 0,$ so that Theorem~\ref{thm1} applies.  We leave
this as an exercise.


In a similar spirit, the function $f(t) = \ln(e^{L}-L+t)$ for 
$t > L - e^L$ is root-like near the fixed point $L$, and so has
properly convergent candidate sequences.  
For an example of a different sort, let $f(t) = 6 - (2t+16)/(t^2+1)$.
This is root-like in a neighborhood of the fixed point at 
$t = 5$, so again candidate sequences for
$f$ converge properly.  

Although Theorem~\ref{thm1} reveals many instances of proper
convergence of candidate sequences, we have not exposed many cases where we can
determine the limits.  For $f(t) = \sqrt{2+t}$ and $t_0 = 0$, an exact
limit of $\pi^2$ occurs as shown in~\eqref{mp3}.  A few more examples
are derived in \cite{geom_ext_ppr}.  In contrast, the proof of
Theorem~\ref{thm1} introduced an infinite family of $Q$ functions for
which exact limits of candidate sequences can be found.  Indeed, the
proof hinges on that as yet unproven fact.  We will see a proof in the
next section, where we will derive properties of our $Q$
functions as part of the larger family of {\em M\"obius} 
functions.

\section{M\"{obius} Functions} \label{mobius sxn}

The ability to evaluate the limit of a candidate sequence for a $Q$
function follows from an explicit expression for the $n$th iterate
$Q^{(n)}$.  This circumstance, which is not typical for the iterates of an
arbitrary function, is observed for the broader class of fractional
linear functions, defined by an equation of the form
\begin{equation} \label{genmobeqn}
 f(t) = \frac{at+b}{ct+d} 
\end{equation}
where $a, b, c, $ and $d$ are real constants.  As mentioned earlier,
with complex constants and a complex variable, \eqref{genmobeqn}
defines M\"{obius} transformations, prompting us to refer to them in
the real context as M\"{obius} functions.

Notice that the function $f$ defined in\eqref{genmobeqn}
reduces to a linear 
function if $c = 0,$ which is not of much interest in this
paper.  So we will assume $c \ne 0$, and without loss take $c = 1$.
Also, we exclude the cases with $a = b = 0$, $a = c = 0,$  $b = d =
0,$ and  $bc = ad$ because they imply that $f(t)$ is constant.  Logically, it is
sufficient to exclude just one condition, $bc = ad,$ which also
excludes the others. Thus, we will be interested in functions of the form
\begin{equation} \label{gen.frac.lin.eq}
f(t) = \frac{at+b}{t+d},
\end{equation}
where $b \ne ad$.  Hereafter, that is what will be meant whenever we
refer to {\em M\"obius functions}. 

As also mentioned earlier, there is an extensive literature on
M\"{obius} transformations.  Some of the results derived below 
concerning fixed points, eigenvalues, and eigenvectors are explicitly
or implicitly already known from the complex context.   However,
rather than citing multiple 
prior works covering a variety of contexts, we will rederive the
specific properties we need.  We do so in the hope that this will
better serve the interests of the readers, and because the derivations are
straight-forward and reasonably brief. We make no claims of originality. 

It is natural to consider M\"{obius} functions  in the context of
generalized mysterious patterns.  As we 
will discuss presently, diagonalizability of the matrix $M_f$ permits
determination of an explicit formula for the iterate $f^{(n)}(t)$.
This in turn permits the evaluation of limits of candidate
sequences.  Even when $M_f$ is not diagonalizable, matrix analysis
leads to limits of candidate sequences, but we will not discuss the
details in this paper.

Note that the $Q$-functions
constitute a proper subset of M\"obius functions.  For example,
the derivative of a $Q$-function is always positive whenever the
function is defined.  That is not true in general of
M\"obius functions.  We will see later that $Q$ must also
have distinct fixed points.

\begin{figure}[hbt]
\centering 
\includegraphics[scale=0.7]{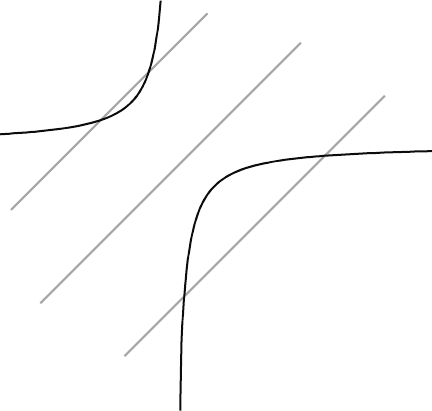}
\caption{When $ad>b$, the graph of $f$ is an equilateral hyperbola,
with an equation of the form 
$(x-h)(y-k) = w$ where $w < 0.$  The gray lines show
possible positions of the line $y=x.$}
\label{mobius.fig1} 
\end{figure}

If we rewrite \eqref{gen.frac.lin.eq} as $f(t) = a + (b-ad)/(t+d)$
we obtain the graph equation $(x-h)(y-k) = w$ where $h = -d,$ $k = a$,
and $w = (b-ad)$.  This is illustrated in \figref{mobius.fig1} for the
case that $w < 0.$  
The fixed points of $f$ are intersections with the
line $y=x$, three possible positions of which are shown in the
figure.  Depending on where $y=x$ sits relative to the curve, there
may be no intersections, a single double intersection, or two
intersections.  In the context of Theorem~\ref{thm1}, where our
approximating M\"obius function is actually a $Q$-function,
there must a fixed point where the curve is concave-down with a slope
in the interval $(0,1)$.  Therefore, the configuration must correspond
to the lower arc in the figure.  In addition, at one fixed point the
curve will have slope less than one, so the intersection with $y=x$
cannot be a tangency.  This implies distinct fixed points, and so the
configuration illustrated by the lowest position of the line $y=x$.
If the fixed points are $L_1 < L_2$, we see that
$f'$ is positive at both and that $f'(L_1)>f'(L_2)$.  These properties
will be useful in understanding iteration sequences 
$\left\{ f^{(n)}(t_0)  \right\}$.

\subsection*{Matrix Analysis Example}
We turn now to an analysis of 
\begin{equation} \label{Qmat}
M_f = \matrrvv{a}{b}{1}{d}
\end{equation}
corresponding to $f(t)$ defined in \eqref{gen.frac.lin.eq}.
If $f$ and $g$ are both M\"{obius} functions, then so is their
composition $f\comp g$ and direct computation shows that $M_{f \circ
g} = M_f M_g$.  This implies 
that the $n$th iterate of $f$ is a M\"{obius} function with
matrix $M_f^n$, providing a means to obtain an explicit
representation of the $n$th iterate of $f$.  

To illustrate, we consider the example 
\[ f(t) = \frac{3t+6}{t+4},  \]
which happens to be a $Q$-function.  Its matrix is
\[ M_f = \matrrvv{3}{6}{1}{4}   \]
and the matrix for $f^{(n)}$ is $M_f^n$.  We compute it using the
diagonalization 
\[ M_f = \frac15 \matrrvv{-3}{2}{1}{1}
                 \matrrvv{1}{0}{0}{6}
                 \matrrvv{-1}{2}{1}{3} = PDP^{-1}.   \]
Since $(PDP^{-1})^n = PD^nP^{-1}$ we have
\begin{eqnarray*}
M_f^n   & = & \frac{1}{5} \matrrvv{-3}{2}{1}{1}
                 \matrrvv{1}{0}{0}{6^n}
                 \matrrvv{-1}{2}{1}{3} \\
   & = & \frac{1}{5}\matrrvv{2\cdot6^n+3}{6(6^n-1)}{6^n-1}{3\cdot6^n+2}.
\end{eqnarray*}
Therefore,
\begin{equation} \label{f^(n)eq1} 
f^{(n)}(t) =
\frac{(2\cdot6^n+3)t+6(6^n-1)}{(6^n-1)t+3\cdot6^n+2}. 
\end{equation}

Using this result, we can proceed to find limits of candidate
sequences.  As a first step, let us find the fixed points of $f$.
They are the solutions to $f(t) = t$, namely $-3$ and 2.  Next we
calculate $f'(-3) = 6$ and $f'(2) = 1/6.$  In particular, $0 < f'(2) < 1$
so $f$ is a contraction on a neighborhood $N$ of 2.  Using our earlier
notations, we have $L = 2$ and $m = 1/6.$

For an initial $t_0 \in N$, the candidate sequence is given by
\[ c_n = \frac{|2-f^{(n)}(t_0)|}{1/6^n} = 6^n|2-f^{(n)}(t_0)|.  \] 
Substituting~\eqref{f^(n)eq1} and simplifying leads to
\[ c_n = \left| \frac{5t_0-10}{t_0+3+(2-t_0)/6^n}  \right|.  \]
Now it is evident that
\[ \lim_{n\rightarrow \infty} c_n = \frac{5t_0-10}{t_0+3}. \]

These results are typical for candidate sequences of $Q$-functions,
and we shall see that the entire analysis can be carried out in 
general.  In particular, by Lemma~\ref{distinct-eigs-lemma}, the
matrix $M_Q$ for any $Q$-function is diagonalizable because it has
distinct real eigenvalues.

\subsection*{Distinct Eigenvalues}
In general there is a connection between fixed points of a M\"obius
function $f$ and the eigenvalues and eigenvectors of $M_f.$
Specifically, $f$ has a fixed point $L$ iff
$L+d$ is an eigenvalue of $M_f$.  To see
this, observe that the fixed point condition for $f$,
\[ \frac{aL+b}{L+d} = L,   \]
implies
\[ aL+b = L(L+d).   \]
Then
\[ \matrrvv{a}{b}{1}{d} \colvv{L}{1} 
=  \colvv{aL+b}{L+d} 
=  \colvv{L(L+d)}{L+d} 
=  (L+d)\colvv{L}{1}. \]
Conversly, suppose $L+d$ is an eigenvalue of $M_f$.  Since by
definition of M\"obius function $M_f$ is non-singular, $L+d \ne 0.$
And $L+d$ must also be a root of
the characteristic polynomial, leading to $L(L+d)=aL+b.$  Thus $f(L)=L.$  

Next consider a $Q$-function $Q$, with distinct eigenvalues, 
$\lambda_1 < \lambda_2$.  It has distinct fixed points $L_1 =
\lambda_1 - d < \lambda_2 - d = L_2.$  At these fixed points we also define $m_j =
Q'(L_j)$, which can be calculated to be $(ad-b)/\lambda_j^2$. This
implies that $m_1m_2 = 1$ because the product of the eigenvalues of
$M_Q$ is equal to the determinant.  We neglect the case that
$|m_1|=|m_2| = 1,$ which is not of immediate interest.  

By direct calculation the fixed points satisfy the equation $aL+b = L(L+d)$ and so 
are given by
\[ L_1 = \frac{a-d - \sqrt{(a-d)^2+4b}}{2} \;\;\; \mbox{and} \;\;\;
   L_2 = \frac{a-d + \sqrt{(a-d)^2+4b}}{2}.   \]
From the definition of $Q$-function we know that at one of the
fixed points the derivative is positive and less than 1.  As shown in
\figref{mobius.fig1}, that must be at the greater of the fixed points,
$L_2$.  So $0 < m_2 < 1$ and it follows that $m_1 = 1/m_2 > 1.$

Now we reverse course, reasoning from the properties of the fixed
points back to the eigenvalues.  We know they are distinct, and given
by $d+L_j$.  Implicit in earlier remarks, we also know that 
corresponding eigenvectors can be expressed in the form $[L_j \;\; 1]^T$.
With this information we can proceed to diagonalize $M_Q,$ and
determine an explicit expression for $Q^{(n)}(t_0).$ 
We formalize the result as the following Lemma.

\begin{lemma} \label{Qfunc orbit}
Let $Q$ be a $Q$-function defined as in~\eqref{Qdef1}.  Then for 
$n\ge 1,$ $Q^{(n)}$ is given by
\begin{equation} \label{Q^{(n)}(t)eqn}
Q^{(n)}(t) =  \frac{b (t - L) m^n + L (b + L t)}
               {L (L - t) m^n+ (b + L t)},
\end{equation}
where $L$ is the greater of $Q's$ two fixed points and $m = Q'(L) \in (0,1)$.
\end{lemma}

Before presenting the proof, we pause to observe that $Q^{(n)}$ is a
M\"obius function, and so defined everywhere except where the
denominator vanishes.  For $n=1$ that occurs at $t = -d$.  For $n=2$
the excluded $t$ satisfies $Q(t) = -d$ so $t = Q^{-1}(-d)$.  And in
general, $Q^{(n)}(t)$ is undefined for 
$t = \left(Q^{-1}\right)^{(n-1)}(-d)$. For future reference we denote
this value of $t$ as $w_n$.  It is the
unique root of the denominator of~\eqref{Q^{(n)}(t)eqn}.

\begin{proof}
To analyze the iterates of $\displaystyle Q(t) = \frac{at+b}{t+d}$ we
will diagonalize
\[ M_Q = \matrrvv{a}{b}{1}{d}.   \]
While we know that $Q$ has fixed points $L_1 < L_2$, we agree to refer to
the greater as $L$, and that $m = Q'(L) \in (0,1)$.  Let us 
also define $s = Q''(L) < 0$.  Applying system~\eqref{mLs to abd}
expresses $a, b,$ and $d$ in terms of $L, m,$ and
$s.$ Thus the characteristic polynomial of $M_Q$ is 
\begin{eqnarray*}
p(t)   & = & t^2 - (a+d)t + ad - b    \\
   & = &  t^2 +\frac{2m}{s}(m+1)t + \frac{4m^3}{s^2}  \\
   & = &  \left(t + \frac{2m^2}{s} \right)  \left(t + \frac{2m}{s} \right).  
\end{eqnarray*}
Evidently the eigenvalues are $\lambda = -2m/s$ and $\mu = -2m^2/s$.
Since $\mu / \lambda = m$, we have $0 < \mu < \lambda.$  This shows
that $\lambda = L+d$ and $\mu = L_1 + d.$

Continuing with the analysis, we have already seen that $[L \;\; 1]^T$
is an eigenvector for $\lambda$, and $[L_1 \;\; 1]^T$ is an
eigenvector for $\mu$. Also, as $L$ and $L_1$ are the roots of
$t^2 + (d-a)t - b$, we see $LL_1 = -b.$ Thus $L_1 = -b/L$ and 
$[b \;\; -L]^T$ is also an eigenvector for $\mu.$

This last conclusion depends in part on our assumption that $L > 0$
for root-like functions.  However if $L = 0,$ which implies $b = 0$ as
well, then $[b \;\; -L]^T$ is the zero vector.  On the other hand,
when this occurs $M_Q$ is triangular with eigenvalues $a$ and $d$,
simplifying the diagonalization process.  We comment that in this
special case it remains possible to determine an explicit expression
for $Q^{(n)}(t)$ as a function of $n$ and $t.$  The details are
left as an exercise. 

Using the expressions derived above for the eigenvalues and
eigenvectors of $M_Q$, we obtain the diagonal representation
\[ M_Q = 
  \frac{1}{b+L^2}\matrrvv{b}{L}{-L}{1}\matrrvv{\mu}{0}{0}{\lambda}
               \matrrvv{1}{-L}{L}{b}.   \]
This implies 
\begin{eqnarray*}
 M_Q^n  
   & = & \frac{1}{b+L^2}\matrrvv{b}{L}{-L}{1}
                        \matrrvv{\mu^n}{0}{0}{\lambda^n}
                        \matrrvv{1}{-L}{L}{b}   \\
   & = &  \frac{1}{b+L^2}\matrrvv{L^2\lambda^n+b\mu^n}{bL(\lambda^n-\mu^n)}
                                 {L(\lambda^n-\mu^n)}{b\lambda^n+L^2\mu^n}
\end{eqnarray*}
Therefore
\begin{eqnarray*}
Q^{(n)}(t)
   & = &  \frac{(L^2\lambda^n+b\mu^n)t +bL(\lambda^n-\mu^n) }
               {L(\lambda^n-\mu^n)t +b\lambda^n+L^2\mu^n }\\
   & = &  \frac{b (t - L) \mu^n + L (b + L t)\lambda^n}
               {L (L - t) \mu^n+ (b + L t)\lambda^n }  \\
   & = &  \frac{b (t - L) (\mu/\lambda)^n + L (b + L t)}
               {L (L - t) (\mu/\lambda)^n+ (b + L t)}  \\
   & = &  \frac{b (t - L) m^n + L (b + L t)}
               {L (L - t) m^n+ (b + L t)}.
\end{eqnarray*}
That is what we wished to show.
\end{proof}

Because $Q$ is root-like, we know that for all $t_0$ sufficiently
close to $L$, the sequence $\left\{ Q^{(n)}(t_0) \right\}$ approaches
limit $L$.  But there are two points worth making.  First, the sequence is not even
defined if $t_0 = w_n$ for some $n$.  Clearly, if $-d < t_0 < L$ then
the iterates $t_n = Q^{(n)}(t_0)$ increase away from $-d$ toward
$L$. Similarly, if $L < t_0 < -d$, then the terms $t_n$ decrease away
from $-d$.  This shows that as long as $t_0$ is closer to $L$ than
$-d$ is, all of the iterates $t_n$ are defined, so $t_0$ is not equal to
$w_n$ for any $n$.  This shows even taking the $w_n$'s into account,
the sequence $\left\{ Q^{(n)}(t_0) \right\}$ converges to $L$ for all
$t_0$ sufficiently close to $L$.

Second, if $t_0$ is either of
the fixed points $L$ and 
$-b/L$, $\left\{ Q^{(n)}(t_0) \right\}$ converges trivially because it
is a constant sequence.  This is correctly indicated
by~\eqref{Q^{(n)}(t)eqn}.  In any other case, as long as all the
iterates $t_n$ exist, ~\eqref{Q^{(n)}(t)eqn}
implies that $\left\{ Q^{(n)}(t_0) \right\}$ converges to $L$.

Now we turn our attention to candidate sequences for $Q$, in the form
of the next Theorem.

\begin{theorem} \label{Qfunc canseq}
Let $Q$ be a $Q$-function defined as in \eqref{Qdef1}, and hence a
contraction on an interval $I$ containing a fixed point $L$ where
$m=Q'(L) \in (0,1)$.  Then for $t_0 \in I$ and $t_0 \ne L$ the candidate sequence
\[ c_n = \frac{|L - Q^{(n)}(t_0)|}{m^n} \]
converges properly, with 
\[ \lim_{n \rightarrow \infty} c_n =
       \left|\frac{(L^2+b)(L-t_0)}
                  {b + L t_0} \right|. \]
\end{theorem}

\begin{proof}
We first use~\eqref{Q^{(n)}(t)eqn} to compute
\begin{eqnarray*}
m^n c_n   
& = & \left|L - \frac{b (t_0 - L) m^n + L (b + L t_0)}
                     {L (L - t_0) m^n+ (b + L t_0)} \right|   \\
& = & \left|\frac{L(L (L - t_0) m^n + (b + L t_0)))
                        -(b (t_0 - L) m^n + L (b + L t_0))}
                     {L (L - t_0) m^n+ (b + L t_0)} \right|   \\
& = & \left|\frac{(L^2+b)(L-t_0)m^n}
                 {L (L - t_0) m^n+ (b + L t_0)} \right|.  
\end{eqnarray*}
Dividing by $m^n$ thus shows
\[ c_n =  \left|\frac{(L^2+b)(L-t_0)}
                     {L (L - t_0) m^n + (b + L t_0)} \right|,  \]
which evidently implies
\[ \lim_{n \rightarrow \infty} c_n  =  
       \left|\frac{(L^2+b)(L-t_0)}
                  {b + L t_0} \right| \]
for $t_0 \ne -b/L$.

Although the limit of $c_n$ is undefined if $t_0
= -b/L$, that is the fixed point of $Q$ where $Q' > 1$, and so not an
element of $I$.  For any other
value of $t_0,$ the candidate sequence converges to a non-negative limit.
That limit is nonzero because $t_0 \ne L$ and $L^2 +b =
2mL(m-1)/s$ where $L \ne 0$ and $m$ is neither 0 nor 1.  Thus, we have
shown that the candidate sequence always converges to a finite
positive limit.
\end{proof}

Using the relations in \eqref{mLs to abd} the limit can be expressed in the form
\[ \lim_{n \rightarrow \infty} c_n  =  
   \left| \frac{1}{t_0-L} + \frac{s}{2m(m-1)} \right|^{-1}. \]
Thus, when $t_0$ is far from $L$, the first fraction above is small,
so that the value of the limit is close to the reciprocal of the
second fraction.  On the other hand, when $t_0$ is very close to $L,$
the first fraction will be very large, and will dominate the second
fraction.  In this case, the limit will be close to zero.

\subsection*{Continued Fractions}
In the world of iterated functions nested squareroots of the sort we
have seen above arise naturally and have an attractive algebraic
appearance.  Another familiar family of iterated functions is
connected with continued fractions.  For example, we might consider
something of the form
\[ 2+\frac{15}{2+\frac{15}{2+\frac{15}{\cdots}}}    \]
Observe that this expression can be obtained by iterating the function 
$f(t) = 2 + 15/t$, and numerically it is clear that the limit is
$L=5$.  This is a fixed point of $f,$ and $f'(L) = -3/5.$  We can
formulate the companion sequence 
\[ c_n = \left(\frac53\right)^n |5 - f^{(n)}(t_0)|.   \]
Can we find the limit?

Yes we can.  Note that $f(t)$ is a M\"{obius} function with 
\[ M_f = \matrrvv{2}{15}{1}{0}   \]  
and although this is not a $Q$-function, the analysis we used above
can be adapted to this case as well.  As before $M_f$ can be
diagonalized, leading to an expression for $c_n$ as a function of $n,$
and the conclusion of Theorem~\ref{Qfunc canseq} still applies.  Thus
\[ \lim_{n \rightarrow \infty} c_n =
       \left|\frac{8(5-t_0)}{3 + t_0} \right|. \]
Taking $t_0 = 2$ we can express the limit in the form of a
mysterious pattern:
\[ \lim_{n \rightarrow \infty}  \left(\frac53\right)^n 
\left|3-\frac{15}{2+\frac{15}{2+\cdots + \frac{15}{2}}} \right| 
= \frac{24}{5}. \] 

This is in a way an analog of the original mysterious pattern.
Similar results can be obtained for any continued fraction with
constant coefficients in place of the values 15 and 2 in the example.

As a particular instance, let us consider the continued fraction that results
from iterating $f(t) = 1+1/t$.  As is well known, extending this
continued fraction for a finite number of stages results in a ratio of
successive Fibonacci numbers.  Following the same analysis as above,
we see that the fixed point is $\phi$, here representing the golden
ratio $(1+\sqrt{5})/2$.  We also compute $m = f'(\phi)$ and find that
$|1/m| = \phi^2$.  In the end we obtain this result:
\[ \lim_{n \rightarrow \infty}  \phi^{2n}
\left|\phi - \frac{F_{n+1}}{F_{n}} \right| = \sqrt{5} \]
where $F_n$ is the $n$th Fibonacci number starting from $F_0 = 0$ and
$F_1 = 1$.

Thus we see that the computability of M\"{obius} function
candidate sequences is not useful solely in proving
Theorem~\ref{thm1}.  It also reveals an entire family of
mysterious-pattern-like results.  And in the case of the second
example above, we obtain a result that doesn't even appear on the
surface to involve function iteration.  It follows the familiar
pattern of an exponential factor multiplied by the difference between
a convergent sequence and its limit.  We only recognize the appearance
of an iterated function when we relate the ratios of successive
Fibonacci numbers to continued fractions.

\subsection*{Dynamic Parameterization of $f$}

We have already seen how to construct a M\"{obius} function $Q$ with
$Q(L) = L,$ $Q'(L)=m,$ and $Q''(L) = s$ using system~\eqref{mLs to
abd}.  Therefore, for an arbitrary function $f$ with a fixed point at
$L$ we can construct a M\"{obius} function $Q$ so that $Q$ and $f$
agree at $L$ through the second derivative.  That is, we can always
construct a second order M\"{obius} approximation to a contraction at
its fixed point.  In this paper we have demonstrated proper convergence of
candidate sequences for root-like functions using a first order
M\"{obius} approximation.  It remains to be seen what use there may be
for second order approximations.

\section{Eigen-Functions} \label{eigfunc sxn}
We turn next to a different approach to simplifying iterated functions.
Recall that 
the original mysterious pattern was formulated in terms of iterates of
the function $f(t)=\sqrt{2+t}$.  This arises in connection with the
cosine half-angle identity, expressed in the form
\[ 2 \cos(t/2) = \sqrt{2 + 2\cos(t)}.  \]
On the right we have $f(2\cos{t})$ and on the left we have
$2\cos(t/2)$.  If we define functions $\phi(t) = 2\cos(t)$ and
$\mu_{1/2}(t) = t/2$, then the half-angle identity becomes
\begin{equation} \label{MP phimu=fphi2}
\phi \comp \mu_{1/2} = f \comp \phi 
\end{equation}
where $\comp$ as usual represents function composition.  This in turn
leads to 
\begin{equation} \label{MP f=phimuphiinv}
f = \phi \comp \mu_{1/2} \comp \phi^{-1}
\end{equation}
(on an interval on which $\phi$ is invertible).  This equation echoes
the diagonalization we saw for M\"{obius} functions.  Indeed, in terms of the
composition operator, we can view $\phi$ as an eigenvector for $f$,
with eigenvalue $\mu_{1/2}$.  In any case, the preceding equation
shows that
\[ f^{(n)} =  \phi \comp \mu_{1/2}^{(n)} \comp \phi^{-1} 
           =  \phi \comp \mu_{1/2^n} \comp \phi^{-1}.   \]
Evaluating this expression at $0 = \phi(\pi/2)$ we obtain
\[ f^{(n)}(0)
  =  \phi \comp \mu_{1/2^n} \comp \phi^{-1} \comp \phi(\pi/2)
  =  2\cos(\pi/2^{n+1}).   \]
Hence,
\begin{eqnarray*}
 2^{n+1}\sqrt{2-f^{(n)}(0))}
    & = & 2^{n+1}\sqrt{2-2\cos(\pi/2^{n+1})} \\
    & = & 2^{n+2}\sin(\pi/2^{n+2}) \\
    & = & \pi \frac{\sin(\pi/2^{n+2})}
                    {\pi/2^{n+2}}.  
\end{eqnarray*}
In the limit this becomes the mysterious pattern.

The preceding analysis suggests an approach to finding other generalizations of the
mysterious pattern for which candidate sequence limits can be
determined.  The idea is to seek other functions $f$ for which an
equation of the form of~\eqref{MP f=phimuphiinv} holds. 
And actually, our earlier analysis of
M\"{obius} functions exactly corresponds to~\eqref{MP f=phimuphiinv}.
In the matrix equation $M_Q = PDP^{-1}$, $P$ and $D$ can be seen as
matrices $M_{\rho}$ and $M_{\delta}$ for M\"{obius} functions
$\rho$ and $\delta$.  Furthermore, $\delta(t) = \alpha t$ where 
$\alpha$ is the ratio of the
diagonal entries of $D$.  Thus the matrix diagonalization represents
the identity $f = \rho \comp \mu_{\alpha} \comp \rho^{-1}$ .

Looking beyond M\"{obius} functions, let us consider what an equation
of the form of~\eqref{MP f=phimuphiinv} implies about candidate
sequences for the corresponding function $f$.  Accordingly, suppose that
\begin{equation} \label{f=phimuphiinv2}
f = \phi \comp \mu_\alpha \comp \phi^{-1}
\end{equation}
where $0 < \alpha < 1$ and $\phi$ is an invertible function mapping its
domain into a neighborhood of $L,$ an attracting fixed point of $f$.
Then at $L$, we have
$\phi \mu_\alpha \phi^{-1}(L) = L$ (where we have suppressed the $\comp$ notation)
so  $\alpha \phi^{-1}(L) =
\phi^{-1}(L)$.  Then unless $\alpha = 1$, in which case $f$ is the
identity, we must have $\phi^{-1}(L) = 0.$  Thus $L = \phi(0)$.  

In general $\phi(0)$ will be a fixed point of $f.$
It will be an attracting fixed point when $|f'(L)| < 1.$
Now if $\theta$ is in the domain of $\phi,$ we can take 
$t=\phi(\theta)$ to express~\eqref{f=phimuphiinv2} as 
$f(\phi(\theta))=\phi(\alpha \theta)$.
Differentiation produces
\[ f'(\phi(\theta))\phi'(\theta) = \alpha \phi'(\alpha \theta) . \]
At $\theta = 0$ this becomes
\[ f'(L)\phi'(0) = \alpha \phi'(0) . \]
When $\phi'(0) \ne 0$ we conclude that $f'(L) = \alpha.$  So since 
$0 < \alpha < 1,$ $L$ will be an attracting fixed point when 
$\phi'(0) \ne 0$ and candidate sequences can be studied.

Proceeding, we know a candidate sequence for $f$ will be of the form
\[ c_n = \frac{|L - f^{(n)}(t_0)|}{m^n}  \]
where $m = f'(L) = \alpha$.  Let $t_0 = \phi(\theta_0)$.  Then 
\[f^{(n)}(t_0) = \phi(\alpha^n\phi^{-1}(\phi(\theta_0))) 
               = \phi(\alpha^n \theta_0). \]
Thus, the candidate sequence is
\[ c_n = \frac{|L - \phi(\alpha^n \theta_0)|}{\alpha^n}
       = |\theta_0| \left| \frac{\phi(0) - \phi(\alpha^n \theta_0)}
                                {0- \alpha^n \theta_0} \right|.  \]
By inspection, 
\begin{equation} \label{c_n lim eq2}
 \lim_{n \rightarrow \infty}c_n = |\theta_0 \phi'(0)|. 
\end{equation}
Therefore candidate sequences converge in this case and the exact
limit can be evaluated.

Interestingly, in the case of the mysterious pattern, $\phi(\theta) = 2
\cos(\theta)$ and $\phi'(0)=0,$ so the above analysis does not apply.
But the analysis can be extended as follows.  Suppose as before that
$f = \phi \mu_\alpha \phi^{-1}$, and also assume that $\phi'(0) =
0,$ but $\phi''(0) \ne 0.$  Starting again from 
$f(\phi(\theta)) = \phi(\alpha \theta)$, 
differentiating twice produces
\[ f''(\phi(\theta))\phi'(\theta)^2 + f'(\phi(\theta))\phi''(\theta)
     = \alpha^2 \phi''(\alpha \theta).  \]
Therefore, for $\theta = 0,$
\[ f'(\phi(0))\phi''(0)
     = \alpha^2 \phi''(0) \]
and
\begin{equation} \label{m=alpha^2 eqn}
 m = f'(L) = \alpha^2.   
\end{equation}

Now the candidate sequence is
\[ c_n = \left| \frac{f^{(n)}(t_0) - L}{m^n} \right|  
       = \left| \frac{\phi(\alpha^n \phi^{-1}(t_0)) - L}
                     {\alpha^{2n}}\right|  
       = \theta^2 \left| \frac{\phi(\alpha^n \theta_0) - \phi(0)}
                              {\alpha^{2n}\theta^2} \right|. \]
This can be rewritten as
\[ c_n = \theta^2 \left| \frac{\phi(\epsilon_n) - \phi(0)}
                              {\epsilon_n^2}  \right| \]
where $\epsilon_n = \alpha^n \theta_0$.  As $n$ goes to infinity,
$\epsilon_n$ goes to 0.  Thus, 
\[ \lim_{n\rightarrow \infty} c_n 
     = \theta^2 \lim_{\epsilon \rightarrow 0} 
          \left| \frac{\phi(\epsilon)-\phi(0)}{\epsilon^2} \right|. \]

Finally, since $\phi'(0) = 0$, the limit on the right equals
$\phi''(0)/2.$   This can be shown in various ways, including by an
application of L'Hospital's rule.  Perhaps a better approach is to
assume that $\phi$ is represented by a power series at zero, which
points immediately to the more general case that the first $k$
derivatives of $\phi$ all vanish at zero.  Be that as it may, for the
case at hand, we have shown that
\begin{equation} \label{phi'(0)=0 eqn1}
 \lim_{n\rightarrow \infty} c_n 
   = \left| \frac{\theta^2 \phi''(0)}{2} \right|
   = \left| \frac{\phi^{-1}(t_0)^2 \phi''(0)}{2} \right|.
\end{equation}

Applying this in the case of the original mysterious pattern, we have 
$\phi(\theta) = 2\cos{\theta}$ and $t_0 = 0$.  Thus $\theta_0 = \pi/2$,
$\phi''(0) = -2,$ and~\eqref{phi'(0)=0 eqn1} becomes
\[  \lim_{n\rightarrow \infty} c_n = \frac{\pi^2}{4}.  \] 
Thus we recover the mysterious pattern as expressed in~\eqref{mp3.1}.

As the results above show, finding appropriate 
instances of~\eqref{f=phimuphiinv2} may well lead us to
generalizations of the mysterious pattern with computable limits.  How
might such instances be found?  We offer the following three approaches.

\subsection{Conjugation.}
To obtain other instances of~\eqref{f=phimuphiinv2}, observe that
one equation of that form leads to others.
For example, with the functions in~\eqref{MP f=phimuphiinv}, we can formulate
\[ \psi \comp f \comp \psi^{-1} = 
\psi \comp \phi \comp \mu_{1/2} \comp \phi^{-1} \comp \psi^{-1}  \]
for suitable invertible $\psi$.  Taking $\psi$ to be the
exponential function $e^t$, this produces the identity
\[ \exp(\sqrt{2+\ln t}) 
    = \eta \comp \mu_{1/2} \comp \eta^{-1}, \]
where $\eta(\theta)= e^{\phi(\theta)}.$
The corresponding analog to the mysterious pattern involves the
candidate sequences of the function $g(t) = \exp(\sqrt{2+\ln t})$ for
$t_0$ near $g$'s fixed point $L = e^2$.  Note that $g'(L) = 1/4 \in
(0,1)$.  However, in this example the eigen-function is not $\phi(\theta)$ but
$\eta(\theta) = e^{2\cos(\theta)}$.  Because $\eta'(0) = 0 \ne \eta''(0)$, we can
apply~\eqref{phi'(0)=0 eqn1} to find
\begin{equation} \label{c_n lim eq1}
\lim_{n \rightarrow \infty} \left| \frac{e^2-g^{(n)}(t_0)}
                                        {1/4^n} \right|
      = e^2\cos^{-1}(\ln (\sqrt{t_0}))^2.
\end{equation}

Although we can in this fashion generate a variety of analogs of the
mysterious pattern, the results seem to be a bit artificial and
contrived.  Moreover, the nested form of
$g^{(n)}$ is essentially the same as that for $f^{(n)}$.  That is
clear because 
\[ g^{(n)}(t_0) = (\psi f \psi^{-1})^{(n)}(t_0) 
   = \psi (f^{(n)} (\psi^{-1}(t_0))).  \]
Consequently the analog of the mysterious pattern produced by $g$,
which involves an expression of the form
\[ c_n = 4^n \left| e^2 - 
   e^{\sqrt{2+\sqrt{2+ \cdots \sqrt{2+\ln{t_0}}}}}    \right|, \]
appears to be an algebraic echo of the original pattern.
In particular, $t_0 = 1$ produces
\[ \lim_{n\rightarrow \infty}
4^n \left| e^2 -  e^{\sqrt{2+\sqrt{2+ \cdots \sqrt{2}}}} \right|
= e^2\frac{\pi^2}{4}.  \]
This can easily be derived from~\eqref{mp3.1}, as we leave it to the
reader to demonstrate.

Nevertheless, it may be worth considering this approach just a bit
further.  In general terms, we can think of the transition
from $f = \phi \comp \mu_{\alpha} \comp \phi^{-1}$ to
$g = \psi \phi \comp \mu_{\alpha} \comp \phi^{-1} \psi^{-1}$ as a
conjugation.  It transforms $f$ into $\psi \comp f \comp \psi^{-1} =
g$.  In the process $\phi$ gets {\em translated} into $\psi \comp
\phi.$  

We can consider how candidate sequences and their limits are affected
under conjugation when $\psi$ has a simple form.  For example, let
$\pi_{\beta}(\theta) = \theta+\beta$.  Then 
$\pi_{\beta}\comp f \comp \pi_{\beta}^{-1}(\theta) = \beta+f(\theta-\beta)$ and
the eigen-function is $\pi_{\beta}\comp \phi(\theta) = \phi(\theta)+\beta$, the
derivative of which is just $\phi'(\theta)$.  If $L$ is a fixed point
of $f$ then $L+\beta$ is the corresponding fixed point of
$\pi_{\beta}\comp f \comp \pi_{\beta}^{-1}$.  And $(\pi_{\beta}\comp f
\comp \pi_{\beta}^{-1})'(L+\beta) = f'(L)$.  Assuming $\phi'(0) \ne 0$, we
can apply~\eqref{c_n lim eq2} to find candidate sequence limits for 
$\pi_{\beta}\comp f \comp \pi_{\beta}^{-1}$:
\[ \lim_{n \rightarrow \infty}c_n = |(\theta_0) \phi'(0)| \] 
where $\theta_0 = \phi^{-1}(t_0-\beta)$.

A similar analysis can be completed for 
$\mu_{\beta}\comp f \comp \mu_{\beta}^{-1}(\theta) = \beta
f(\theta/\beta)$.  Without repeating all the details, we note that the
candidate sequence limit in this case turns out to be
\[   |\theta_0 \beta \phi'(0)| \] 
where $\theta_0 = \phi^{-1}(t_0/\beta)$.

\subsection{$\phi$ First.}
Our second method for obtaining additional instances
of~\eqref{f=phimuphiinv2} is to choose first the eigen-function
$\phi$ and then find a corresponding function $f.$  This will succeed
for a given 
$\phi$ if $\phi(\alpha \theta)$ can be expressed in terms of
$\phi(\theta)$, represented by the function $f.$  That implies
$\phi(\alpha \theta) = f(\phi(\theta))$, an alternative form of~\eqref{f=phimuphiinv2}.

As an example, define $\phi(\theta) = e^{\theta}.$ Then since
\[ \phi(\alpha \theta)) = e^{\alpha \theta} = \left( e^{\theta}
\right)^{\alpha} = \phi(\theta)^{\alpha},   \]
we see that $\phi(\alpha \theta) = f(\phi(\theta))$ holds with 
$f(t) = t^{\alpha}$.  Now this $f$ has a fixed point $L$ at $t=1$
and $f'(L) = \alpha.$  When $0 < \alpha < 1$, $f$ is a contraction
in a neighborhood $N$ of $L$ and we can formulate candidate sequences for $f$.  In
addition, $\phi'(0) = e^0 = 1,$ so for this example~\eqref{c_n lim
eq2} does apply.  It shows that for any $t_0 \in N$ the candidate sequence must converge
to $|\phi^{-1}(t_0)\phi'(0)| = |\ln(t_0)|$.

Further reflection shows that for this example we have $f^{(n)}(t_0)
= t_0^{\alpha^n}$.  Thus we can express $c_n$ as an explicit function
of $n$ and calculate the limit directly.  Although this does show that
the eigen-function method was superfluous here, it is at least a
comfort to demonstrate that it produces the correct results in this
case.

For the preceding example we were free to choose $\alpha$ anywhere in
the interval $(0,1)$.  In other examples the function $\phi$ might
only satisfy $\phi(\alpha \theta) = f(\phi(\theta))$ for specific
values of $\alpha.$  Indeed we have already observed that when
$\phi(\theta) = 2 \cos (\theta)$ it is possible to express
$\phi(\alpha \theta)$ as a function of
$\phi(\theta)$ when $\alpha = 1/2$.  Likewise for any integer
$\alpha$ it is well known that $\cos(\alpha \theta)$ can be expressed in terms of
$\cos(\theta)$ using Chebyshev polynomials.  This fact leads to
additional extensions of the mysterious pattern, as we shall see in
Section~\ref{Cheby_stuff}. 

As a variant of the $\phi(\theta) = e^{\theta}$ example, we suggest
defining $\phi(\theta) =\frac12 e^{\theta}$ and $\alpha
= 1/3$.  The reader is invited to verify that:
\begin{itemize}
\item $\phi(\alpha \theta) = f(\phi(\theta))$ holds for $f(t) = \sqrt[3]{(1/4)t}$.
\item $f$ has a fixed point at $L = 1/2.$
\item $f'(L) = 1/3.$
\item $\phi'(0) \ne 0$
\item The candidate sequence associated with $f^{(n)}(t_0)$ converges to 
$|\ln(2t_0)|/2$ for $t_0$ near $L$.
\end{itemize} 
The final assertion can be shown either by using~\eqref{c_n lim eq2} or by
expressing $c_n$ as an explicit function of $n$ and computing the
limit directly.  In either case, the value of the limit can be
expressed in the form
\[ \lim_{n \rightarrow \infty} 
  3^n \; \rule[-16pt]{.5pt}{48pt} \;  \frac12 - \sqrt[3]{\frac14\times\sqrt[3]{\frac14\times
  \cdots \times \sqrt[3]{\frac14}}}\;\rule[-16pt]{.5pt}{48pt}
  = \frac12 |\ln(2t_0)|    \]
where $n$ is the number of nested radicals.  This is clearly an analog
of the original mysterious pattern.

\subsection{Constructing Eigen-Functions.} \label{find phi sxn}
Our third (and possibly favorite) way to obtain instances
of~\eqref{f=phimuphiinv2} is to
begin with some function $f$ having properly convergent candidate sequences
and attempt to construct an appropriate eigen-function $\phi$.  In
particular, we try to produce a power series for the desired function
$\phi.$

Before proceeding, we observe that $\phi$ need not be uniquely
determined by $f$.  A sufficient condition for a function $f$ to have two decompositions
\[f = \phi_1 \mu_{\alpha} \phi_1^{-1} 
    = \phi_2 \mu_{\alpha} \phi_2^{-1}  \]
is for $\phi_1(t) = \phi_2(\beta t)$ for some constant $\beta$.  If
so, then $\phi_1 = \phi_2 \comp \mu_{\beta}$ so
\[  \phi_1 \mu_{\alpha} \phi_1^{-1}
   = \phi_2 \mu_{\beta}\mu_{\alpha} \mu_{\beta}^{-1}\phi_2^{-1}
   = \phi_2 \mu_{\alpha} \phi_2^{-1}.  \]
This holds because 
$\mu_{\beta}\mu_{\alpha} = \mu_{\beta \alpha} = \mu_{\alpha}\mu_{\beta}$.  
Consequently, if we can express $f$ in terms of a function $\phi_1$
with $\phi_1'(0)\ne 0$, then we can also express $f$ in terms of a function
$\phi_2$ with $\phi_2'(0) = 1$.  Thus we need only consider two cases
in attempting to construct $\phi$: either $\phi'(0) = 0$ or $\phi'(0)
= 1$.  In either case we can assume that $\phi$ is analytic in a
neighborhood of $f'$s fixed point $L$ and try to construct its power
series there.  The construction we use is essentially the same as
Devaney's~\cite{devaney}, though we have a different objective.

To illustrate the construction, we consider the example $f(t) = \sqrt{6+t}$ with $L$ = 3.
Expressing~\eqref{f=phimuphiinv2} in the form 
$f \comp \phi = \phi \comp \mu_{\alpha} $ we derive
$\sqrt{6+\phi(\theta)} = \phi(\alpha \theta)$.  
Replacing $\theta$ with $\theta/\alpha$ we find
$\sqrt{6+\phi(\theta/\alpha)} = \phi(\theta)$, or more simply
\begin{equation} \label{L3D0eqn}
6+\phi(\theta/\alpha) = \phi(\theta)^2. 
\end{equation}
This identity can be
differentiated repeatedly to find a power series expansion about 0 for
$\phi$, using the fact that $\phi(0) = L.$  

To begin the process, differentiate~\eqref{L3D0eqn} to obtain
\[ \frac{1}{\alpha} \phi'\left(\frac{\theta}{\alpha}\right) 
      = 2\phi(\theta)\phi'(\theta).   \]
At $\theta=0$ this becomes
\[ \frac{1}{\alpha} \phi'(0)
      = 2\phi(0)\phi'(0) = 2L\phi'(0)   \]
hence either $\phi'(0) = 0$ or else $\phi'(0) \ne 0$ and $1/\alpha = 2L = 6$.
We proceed with the second case, so that~\eqref{c_n lim eq2}
can be applied.  And as we have just shown, we can assume
$\phi'(0) = 1$ without loss of generality.
Accordingly, we set $\phi'(0) = 1,$ and continue with
our program.

Using the fact that $\alpha = 1/6$ we restate~\eqref{L3D0eqn} as
\begin{equation} \label{L3D0eqn2}
6+\phi(6\theta) = \phi(\theta)^2
\end{equation}
and differentiate to obtain
\[ 6\phi'(6\theta) = 2\phi(\theta)\phi'(\theta). \]
Differentiating again produces
\[ 36\phi''(6\theta) =
2\phi(\theta)\phi''(\theta)+2\phi'(\theta)^2. \]
Now set $\theta = 0$ and use our earlier results $\phi(0) = 3$ and
$\phi'(0) = 1$ to find
\[ 36\phi''(0) - 6\phi''(0) = 2. \]
Thus $\phi''(0) = 1/15$ and $\phi''(0)/2 = 1/30$.  In particular, the
second order Taylor polynomial for $\phi$ is 
$p_2(\theta) = 3 + \theta +\theta^2/30.$
We are well on our way!

Succeeding derivatives of~\eqref{L3D0eqn2} become increasingly
complicated.  To help organize the calculations we apply Leibniz's
rule for derivatives of products to $\phi^2 = \phi \cdot \phi.$  It
says that
\[ (\phi\phi)^{[n]}(0) = \sum_{k=0}^n \binom{n}{k} \phi^{[k]}(0)\phi^{[n-k]}(0) \]
where the $[n]$ exponent represents $n$th derivative.  This can be
compressed using symmetry, depending on the parity of $n$.  Thus
\[ (\phi\phi)^{[2n]}(0) = 2\sum_{k=0}^{n-1} \binom{2n}{k}
\phi^{[k]}(0)\phi^{[2n-k]}(0)  + \binom{2n}{n}\phi^{[n]}(0)^2 \]
and
\[ (\phi\phi)^{[2n+1]}(0) = 2\sum_{k=0}^{n} \binom{2n+1}{k}
\phi^{[k]}(0)\phi^{[2n+1-k]}(0). \]
In either case the first term for the $n$th
derivative is $2\phi(0)\cdot\phi^{[n]}(0) = 6\phi^{[n]}(0).$

A relatively simple equation for determining the
value of $\phi^{[n]}(0)$ can be derived using these results.  They
show that the $n$th derivative of~\eqref{L3D0eqn2} at $\theta=0$ is
\[6^n \phi^{[n]}(0) = 6\phi^{[n]}(0) + \sum \mbox{higher order terms}. \]
Thus we find 
\[\phi^{[n]}(0) = \frac{1}{6^n-6} \cdot \sum \mbox{higher order terms}. \]
Two examples will show the pattern of the higher order terms clearly.
\[\phi^{[5]}(0) = 2\frac{\binom{5}{1}\phi^{[1]}(0)\phi^{[4]}(0) +
                         \binom{5}{2}\phi^{[2]}(0)\phi^{[3]}(0)}
                        {6^5-6} . \]
\[\phi^{[6]}(0) = \frac{2\left(\binom{6}{1}\phi^{[1]}(0)\phi^{[5]}(0) +
                         \binom{6}{2}\phi^{[2]}(0)\phi^{[4]}(0)\right)
                        + \binom{6}{3}\phi^{[3]}(0)^2 }
                        {6^6-6} . \]

As these equations illustrate, we can compute the value of each
$\phi^{[n]}(0)$ using previously computed values of $\phi^{[k]}(0)$
with $k<n$.  Thus it is possible to recursively compute successive
values of $\phi^{[n]}(0)$ as far as we wish, and hence find the
coefficients $\phi^{[n]}(0)/n!$ of the power series expansion of
$\phi$.

Having already determined the first three terms of the series, we note
that 
\[\phi^{[3]}(0)/3!
     = 2\frac{\binom{3}{1}\phi^{[1]}(0)\phi^{[2]}(0)}{6(6^3-6)} 
     = \frac{6\cdot 1 \cdot 1/15}{6^2\cdot 35} 
     = \frac{1}{3150}. \]
So the fourth term is $\theta^3/3150$.  Computing two additional terms
gives us the fifth order Taylor polynomial
\[p_{5}(\theta)
  =3+\theta+\frac{\theta^{2}}{30}+\frac{\theta^{3}}{3150} 
    +\frac{11\theta^{4}}{8127000}+\frac{97\theta^{5}}{31573395000}. \]

Can we infer anything from these results?  A pattern in the
coefficients might lead us to an explicit representation of $\phi$
in terms of elementary (and other known) functions.  So far, we have
not detected such a pattern.

Examining the sequence $\{\phi^{[n]}(0)\} = 3, 1, 1/15, 1/525,
11/338625, \ldots$ we see that the values are rapidly decreasing.  In
particular, after the initial 3 all of these values are less than or
equal to 1.  If this pattern continues for all $n$, the series for
$\phi$ will be dominated by $2+e^{\theta}$, and thus will converge for
all $\theta.$  At least this would show that our methods do in fact
define $\phi$ as a smooth function.

Accordingly, we proceed to show that $\phi^{[n]}(0) \le 1$ for 
$n \ge 1$.  Indeed, we already know that this is true for $n = 1$.
Applying Leibniz's rule directly for $n>1$, we can express
\[ \phi^{[n]}(0) 
   = \frac{\sum_{k=1}^{n-1} \binom{n}{k}\phi^{[k]}(0)\phi^{[n-k]}(0)} 
          {6^n-6}.\]
If each factor of the form $\phi^{[j]}(0)$ in the summation is less
than or equal to 1, then the summation is no more than $2^n - 2$.
Thus
\[ \phi^{[n]}(0) 
   \le \frac{2(2^{n-1}-1)}{6(6^{n-1}-1)} < 1,\]
as desired.

There remains some ambiguity about whether the function $\phi$ defined
by the series constructed above actually satisfies the desired
identity~\eqref{L3D0eqn}.  We have defined the series in such a way
that an infinite number of necessary conditions have been met.  But
do these conditions together constitute a sufficient condition?  If we
already knew of the existence of a function satisfying~\eqref{L3D0eqn}
and represented by its Taylor series at zero, then we could argue that
it must have the same series as our constructed $\phi$.  But can we
infer existence of a solution from the existence of the constructed
function?  For us this question remains unanswered.

We have considered the matter numerically, however.  We know that our
constructed $\phi$ is approximated by its Taylor polynomials $p_n$.  So we
may ask whether $p_5,$ for example, nearly satisfies
$f \comp \phi = \phi \comp \mu_{\alpha}$ with $f(t) = \sqrt{6+t}$ and
$\alpha = 1/6$.  That is, does the equation
\[ \sqrt{6+p_5(\theta)} = p_5(\theta/6) \]
hold approximately?

Numerical results are encouraging.  Considering $\theta \in [-2,2]$ we
computed values of the error $E_n(\theta) = |\sqrt{6+p_n(\theta)} -
p_n(\theta/6)|$ for $n = 3, 4, 5.$  We observed that the errors
decreased with $n$, with $E_3 < 4.5\cdot 10^{-6}$, 
$E_4 < 2.1 \cdot 10^{-8}$, and $E_5 < 5.1\cdot 10^{-11}$.  We also
numerically computed the averages
\[ \overline{E}_n = \frac14 \int_{-2}^{2} E_n(t)dt ,\]
finding $\overline{E}_3 \approx 7\cdot 10^{-7}$, $\overline{E}_4 \approx 3\cdot 10^{-9}$, and
$\overline{E}_5 \approx 6\cdot 10^{-12}$.

The numerical results also suggest that for $\theta \in [-2,2]$,
$\phi(\theta)$ is increasing, and hence invertible, with range
approximately  $I = [1.1,5.1]$.  On that interval $f$ is a contraction, and
the sequence $\{f^{(n)}(t_0)\}$ increases to $3$ when $t_0 < 3$ and
decreases to $3$ when $t_0 > 3$.  Now for $\theta_0 \in [-2,2]$, if
we take $t_0 = \phi(\theta_0)$, the identity
$f^{(n)}(t_0) = \phi(\theta_0/6^n)$ will hold.  Thus~\eqref{c_n
lim eq2} shows that the candidate sequence for $f$ will have limit
$|\theta_0 \phi'(0)|$.  This equals $|\theta_0|$ because $ \phi'(0) =
1.$

Suppose $t_0 = 0.$  Then $\phi(\theta_0) = 0,$ but that does not provide a
direct means to find $\theta_0$.  By numerical methods we find
$p_5(\theta) = 0$ for $\theta = -3.3656575319\ldots$.  Working
directly from the definition of the candidate sequence, we numerically
estimate the limit to be $3.3656575397\ldots$.  Our estimate based on
the eigen-function approach is correct to seven decimal places.

\subsection{Original Mysterious Pattern.} \label{origMP sxn}
To gain additional insight about the methods applied in the prior
example, let us consider again the function $f(t) = \sqrt{2+t}$ that
arises in conjunction with the original mysterious pattern.  As
discussed at the beginning of this section, near $f$'s fixed point
$L=2$, the function 
$\phi(\theta)=2\cos(\theta)$ is an eigen-function with corresponding
$\alpha = 1/2$ (see~\eqref{MP phimu=fphi2} 
and~\eqref{MP f=phimuphiinv}).  However, as $\phi'(0) = 0$ in this
example, we were unable to apply~\eqref{c_n lim eq2}.  Observe also
that the condition $\alpha = f'(L)$, valid when $\phi'(0) \ne 0$,
fails to hold here.  Instead, consistent with~\eqref{m=alpha^2 eqn},
we have $f'(L) = 1/(2f(L)) = 1/(2L) = 1/4 = \alpha^2.$

What happens if we try to construct an eigen-function for this $f$
using the Taylor series method applied above?  As before, we start
with the eigen-function equation in the form
\[ f(\phi(\theta)) = \phi(\alpha \theta)  \]
and derive a convenient functional equation, in this case
\begin{equation} \label{L2D0eqn}
2+\phi\left(\frac{\theta}{\alpha}\right) = \phi(\theta)^2. 
\end{equation}
Also as before, we now have two cases to consider, depending on
whether or not $\phi'(0)$ vanishes.

Case 1.  Consistent with the eigen-function $\phi(\theta) = 2 \cos \theta$ discussed
above, we can assume $\phi'(0) = 0$.  Then, continuing to use
foreknowledge, if we take $\alpha = 1/2,$ and $\phi''(0) = -2$, we can
recursively find terms of the power series as before.  The resulting
series,   
\[ \phi(\theta) = 2 + 0 - \frac{2\theta^2}{2!}+ 0 +
\frac{2\theta^4}{4!} + 0 - \frac{2\theta^6}{6!} + \cdots,   \]
does indeed reconstruct our known function $\phi$.  This is hardly
surprising, and in a sense validates the power series approach.
Nevertheless, we proceed to investigate the second case.

\begin{figure}[hbtp]
\centering 
\includegraphics[width=3in]{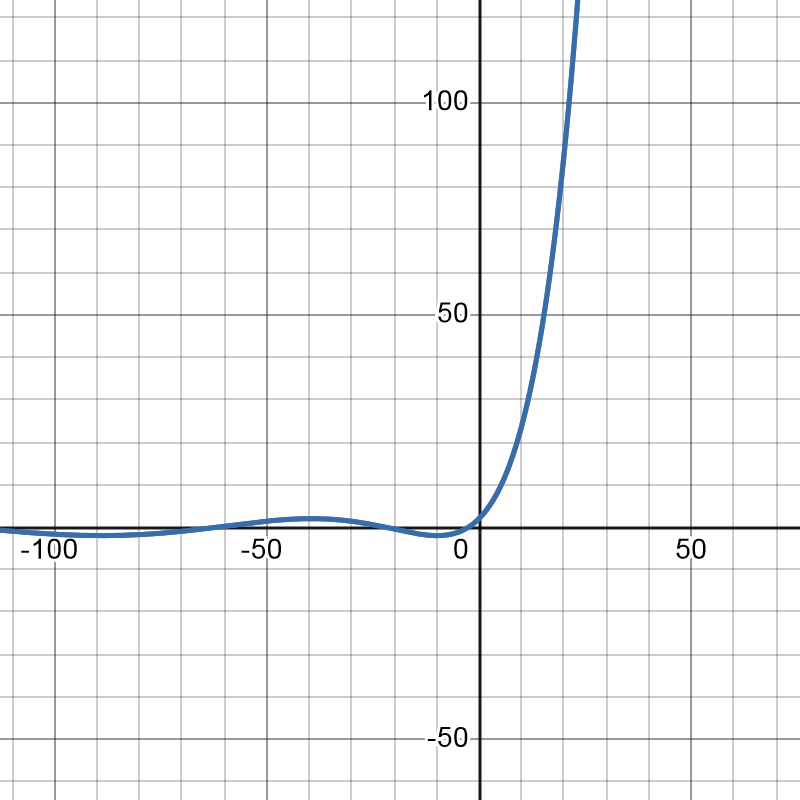}
\caption{Graph of the eigen-function $\phi$ generated as a power
series.  This zoomed-out view shows the hyperbolic cosine like curve
for positive $\theta$ and the oscillating cosine like curve for
negative $\theta.$}
\label{zoomed out graph}
\end{figure}

\begin{figure}[hbtp]
\centering 
\includegraphics[width=3in]{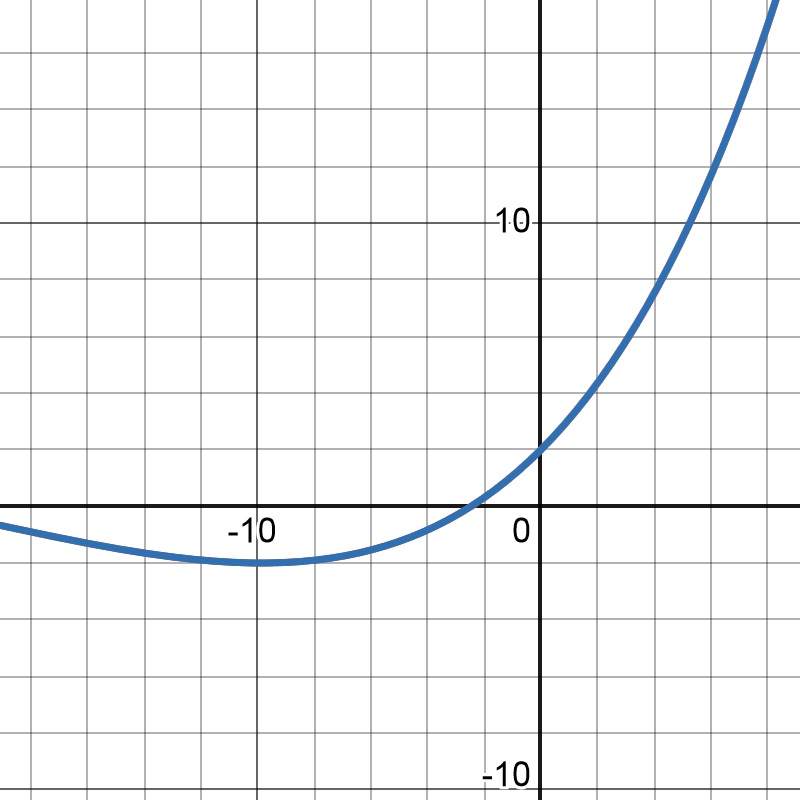}
\caption{Zoomed in graph of $\phi(\theta)$ for $|\theta| < 10$, roughly.  For our
purposes we only need to consider $\theta \ge -\pi^2$.  On that
 interval $\phi$ is increasing.}
\label{zoomed in graph} 
\end{figure}

Case 2.  Let us assume
$\phi(0)\ne 0$.  Then we can proceed with $\alpha = f'(L) = 1/4$, and
assuming $\phi'(0) = 1$, use identity~\eqref{L2D0eqn} and its
derivatives to compute power series coefficients.  The series that
emerges is
\[ \phi(\theta)=2+\theta+\frac{\theta^{2}}{12}+\frac{\theta^{3}}{360}+
\frac{\theta^{4}}{24\cdot840} + \cdots . \]
And, in contrast to the example with $f(t) = \sqrt{6+t}$, this time we
see a pattern: for $k>0$ the $k$th term is $\theta^k/(k(2k+1)!)$.
Assuming that pattern holds for all terms, the expansion
\[ \phi(\theta)=2+\sum_{k=1}^{\infty}\frac{\theta^{k}}{k (2k-1)!}  \]
can be analyzed using standard series methods (or calculus-savvy
software) to find
\[ \phi(\theta) = \left\{ 
   \begin{array}{cl}     
     2\cos\left(\sqrt{-\theta}\right) & \mbox{for} \;\theta<0 \\
     2\cosh\left(\sqrt{\theta}\right) & \mbox{for} \;\theta>0 
   \end{array} \right.   \]
Figures~\ref{zoomed out graph}, \ref{zoomed in graph}, and \ref{phi
with phi inverse}, below, shed light on the nature of $\phi$.

\begin{figure}[hbtp]
\centering 
\includegraphics[width=3in]{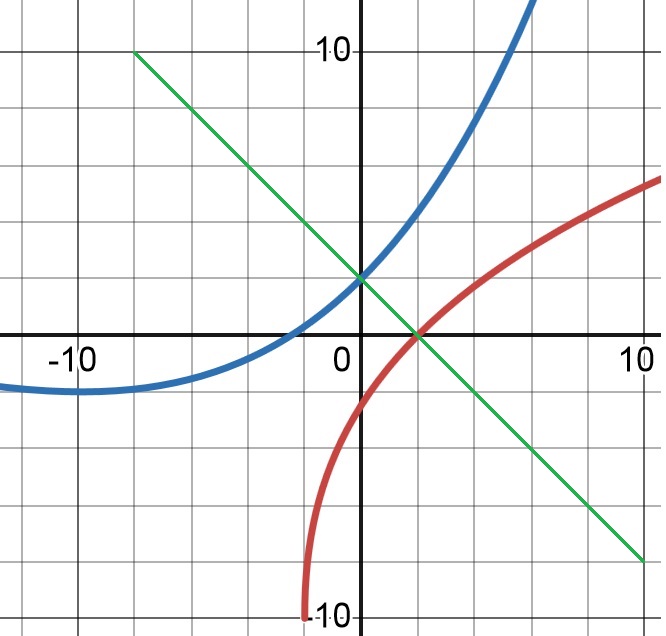}
\caption{Graphs of $\phi$ (in blue) and $\phi^{-1}$ (in red) near
 $(0,0)$.  The line through $(0,2)$ and $(2,0)$ divides the plane into
 two regions.  Below the line $\phi$ and $\phi^{-1}$ are defined in
 terms of cosine; above the line in terms of hyperbolic cosine.}
\label{phi with phi inverse} 
\end{figure}

It turns out that for our purposes, the domain of $\phi$ can be
restricted to $[-\pi^2,\infty)$, on which $\phi$ is increasing.  Then 
there is a corresponding piecewise discription of the inverse function 
\[ \phi^{-1}(t) = \left\{ 
   \begin{array}{cl}     
     -(\cos^{-1}(t/2))^2 & \mbox{for} \;-2 \le t \le 2 \\
      (\cosh^{-1}(t/2))^2 &  \mbox{for} \;t \ge 2 
   \end{array} \right.   \]

Thus we have constructed a second eigen-function for $f(t) =
\sqrt{2+t}$.  And because this new eigen-function has $\phi'(0)=1$, we
immediately obtain the limits of candidate sequences by 
applying~\eqref{c_n lim eq2}.  This produces

\[ \lim_{n\rightarrow \infty} c_n = \left\{ 
   \begin{array}{cl}     
     \cos^{-1}(t_0/2)^2 & \mbox{for} \;-2 \le t_0 \le 2 \\
      (\cosh^{-1}(t_0/2))^2 &  \mbox{for} \;t_0 \ge 2 
   \end{array} \right.   \]

In particular, with $t_0 = 0$ we have a limit of $\pi^2/4$.  This is
consistent with the original mysterious pattern where we computed the
limit not of the candidate sequence, but of a multiple of the
squareroot.  That is, we found
\[ \lim_{n \rightarrow \infty} 2^n \sqrt{2 - f^{(n-1)}(0)} = \pi \]
with $n$ equal to the number of squareroots.  By a trivial modification,
\[ \lim_{n \rightarrow \infty} 2^{n-1} \sqrt{2 - f^{(n-1)}(0)} = \pi/2 \]
which of course is equivalent to 
\[ \lim_{n \rightarrow \infty} 2^{n} \sqrt{2 - f^{(n)}(0)} = \pi/2. \]
After squaring,
\[ \lim_{n \rightarrow \infty} 4^{n} (2 - f^{(n)}(0)) = \pi^2/4. \]
This is none other than the limit of the candidate sequence $c_n$
computed above.

The two $\phi$'s are evidently related.  But we can make the connection
more transparent.  For clarity, let us refer to the original
eigen-function as $\phi_1$ and the more recent eigen-function,
discovered by way of its power series, as $\phi_2$.  Also, for this
discussion we will restrict $t_0$ to the interval $[-2,2)$.  Now
observe that 
\[ \phi_2 = \phi_1 \comp \eta \]
where $\eta(\theta) = \sqrt{-\theta}$ for $\theta \in [-\pi^2,0]$.  Thus we find
\begin{equation} \label{2 phis eq 1}
 \phi_2 \comp \mu_{1/4} \comp \phi_2^{-1}  
   = \phi_1 \comp \eta \comp \mu_{1/4} \comp \eta^{-1} \comp
   \phi_1^{-1}. 
\end{equation}

In the center, the composition $\eta \comp \mu_{1/4} \comp \eta^{-1}$
simplifies to $\mu_{1/2}$ .  This is because $\eta^{-1}(t) = -t^2$
for $t \in [0,\pi]$. So
\[ \eta \comp \mu_{1/4} \comp \eta^{-1}(t)
       = \eta((1/4)\eta^{-1}(t))
       = \eta((1/4)(-t^2))
       = (1/2)|t| = (1/2)t. \]
Hence~\eqref{2 phis eq 1} becomes
\[ \phi_2 \comp \mu_{1/4} \comp \phi_2^{-1}  
   = \phi_1 \comp \mu_{1/2} \comp \phi_1^{-1} = f.     \]
In this way the eigen-equation for $f$ and $\phi_1$ implies the
corresponding equation for $f$ and $\phi_2.$

In a sense, the analysis that produced $\phi_2$ constitutes one instance of
applying the eigen-function method to obtain the limit of candidate
sequences for a given function $f.$  And it is a spectacular success.
We conjure an eigen-function $\phi$ out of thin air, as it were, and then
instantly obtain the candidate sequence limits for all valid $t_0.$
The mysterious pattern limit drops out effortlessly as just one of
an infinite family of such limits.

But in another sense, it is disappointing that this example only
serves to reproduce previously known results.  Even the formulation of
a candidate sequence limit as a function of $t_0$ was already obtained
by geometric methods in~\cite{geom_ext_ppr}.  It would be more
satisfying to produce an example where the eigen-function method
provides a previously unknown candidate sequence limit for some
function $f$.  We will see some examples of this sort later.  But
first we discuss a connection between the methods presented in this
section and a classical result from complex analysis.

\section{Koenigs Functions} \label{Koenigs_stuff}

In the preceding section we described $\phi$ as an eigen-function for
$f$ (under composition) when $f \comp \phi = \phi \comp \mu_{\alpha}$ for some
$\alpha$, or equivalently $f = \phi \comp \mu_{\alpha} \comp
\phi^{-1}$.  This is closely related to an analogous situation in
complex analysis, where $\phi^{-1}$ would be referred to as a {\em 
Koenigs} function for $f$.  So the existence of a Koenigs function
implies the existence of an eigen-function and hence to the
possibility of evaluating limits of candidate sequences.  In this
section we will discuss Koenigs functions and their connection to our
interest in candidate sequences.  It should be emphasized though that Koenigs
functions were developed as a tool for analyzing complex dynamical
systems.  That is distinct from our goal of finding exact limits of
candidate sequences.

Before proceeding further, we pause for an observation inspired by the
Koenigs function material: the existence of a Koenigs function is
implied by proper convergence of candidate sequences.  To understand
this idea, consider a function $f:I \rightarrow I$ for some interval
$I$, such that $f$ has a fixed point $L \in I$, $f'(L)=m,$ and
the candidate sequence $|f^{(n)}(t_0) - L|/m^n$ converges
for all $t_0\in I$.  Then we can define a function $h(t)$ as the
pointwise limit of the sequence of functions 
$h_n(t) = |f^{(n)}(t) - L|/m^n$. Now observe that
\begin{equation} \label{hof=mh}
 h(f(t)) 
= \lim_{n \rightarrow \infty} \frac{|f^{(n)}(f(t))-L|}{m^n}  
= m \lim_{n \rightarrow \infty} \frac{|f^{(n+1)}(t)-L|}{m^{n+1}} 
= m h(t). 
\end{equation}
Thus, where $h^{-1}$ is defined, 
\begin{equation} \label{f=h^{-1}mh}
f = h^{-1} \comp \mu_m \comp h 
\end{equation}
and $h^{-1}$ is an eigen-function for $f$.  

In other words, loosely speaking, not only can we use an eigen-function
to find limits of candidate sequences, the very existence of those
limits implies the existence of an eigen-function. 

These ideas connect with an area of complex analysis concerned with
dynamical systems.  In that setting, for an analytic function $f$
defined on a complex domain, orbits under iteration of $f$ are of
interest.  The functions can be considered in equivalence classes with
respect to conjugation, because under suitable assumptions, the dynamics of $f$ are locally
the same as those of $h\comp f\comp h^{-1}$.  So, for example, one can
always translate a neighborhood of a fixed point $L$ to a neighborhood
of the origin, because translation and scaling are achieved by
conjugation with linear functions.

For the specific situation we have
been studying, where $f$ is a contraction to a fixed point $L$ and 
$|f'(L)|<1,$ the analogous complex case is an analytic function on a disk
centered at $L$.  Then an equation of the form of~\eqref{f=h^{-1}mh} says that
near $L$ $f(t)$ is conjugate to the linear map $\mu_m(t)=mt$, and so
will have the same dynamics. This brings us to a classical result from dynamical systems
that specifies conditions under which~\eqref{f=h^{-1}mh} holds.  It is
due to Gabriel Koenigs, and the function $h$ whose existence it
asserts is referred to as a Koenigs function.

There are various versions of the theorem asserting the
existence of a Koenigs function (see for example \cite[Theorem
2.1]{carleson}), \cite[Theorem 4.1]{devaney}, \cite[Theorem
8.2]{milnor}).  Here we will use the version proven 
in~\cite{carleson}, as stated below.  

\begin{theorem} \label{koenigs-thm}
Suppose $f$ has an attracting fixed point at $z_0,$ with multiplier
$\lambda$ satisfying $0 < |\lambda| < 1.$  Then there is a conformal
map $\zeta = \phi(z)$ of a neighborhood $z_0$ onto a neighborhood of
0 which conjugates $f(z)$ to the linear function $g(\zeta)=\lambda
\zeta$.  The conjugating function is unique, up to multiplication by a
nonzero scale factor.
\end{theorem}
In the original context, it is understood that $f$ is an analytic
function, that the {\em multiplier} is defined to be $f'(z_0)$, and
that a fixed point $z_0$ is said to be {\em attracting} when the modulus of the
multiplier is less than 1.

Beginning by assuming $z_0 = 0$, the proof presented
in~\cite{carleson} has two parts.  One part shows, with $\phi_n(z)$ defined to
be $\lambda^{-n}f^{(n)}(z)$, that the sequence $\{\phi_n(z)\}$ converges
uniformly on a neighborhood of 0. (Observe that $\phi_n(z) =
(f^{(n)}(z)-z_0)/\lambda^n$.) In the other part, defining $\phi(z)$ to
be the 
limit of $\phi_n(z)$, the equation $\phi \comp f = \lambda \phi$ is
derived.  The argument is essentially what appears in~\eqref{hof=mh}.
The uniqueness claim is established by reference to an earlier
result.  There is also an additional comment indicating that the
construction in the proof implies $\phi'(0) = 1.$ 

Now we see that Theorem~\ref{koenigs-thm} may be applicable
to the context of real functions and candidate sequences.  The idea is
to extend $f$ to a complex domain.  We will discuss this in greater detail presently.  But
supposing that $f$ can be extended to an analytic complex function,
the proof of Theorem~\ref{koenigs-thm} immediately implies that
candidate sequences converge properly for $t_0$ in some neighborhood
of the fixed point $L$.  Then, applying Lemma~\ref{lemma3}, the same
conclusion holds for all $t_0$ in the interval containing $L$ on which
$f$ is a contraction.  Indeed, it must hold throughout the basin of
attraction for the fixed point $L$.

Although this provides a possible alternative to Theorem~\ref{thm1}
for establishing proper convergence of candidate sequences, it does
not necessarily help us evaluate their limits.  It is true that such a
limit must be given as $\phi(t_0)$ (which is actually $\phi^{-1}(t_0)$
according to our earlier meaning of $\phi$), but that is of little
value unless we can actually compute $\phi(t)$.  In the context of  
Section~\ref{find phi sxn}, applying Theorem~\ref{koenigs-thm} may tell
us that an eigen-function exists.  But to evaluate the eigen-function
we must be able to invert a function whose values come from the very family of limits we
are hoping to find.  The logic is unfortunately circular, or even
nontrivially knotted.

As given by the proof outlined above Theorem~\ref{koenigs-thm} is an
existence result.  At its heart it applies the Cauchy criterion to
infer the existence of the pointwise limit of a sequence of
functions.  So in general the theorem by itself will not lead to an
explicit representation of the function $\phi.$

On the other hand, in~\cite{devaney} a version of
Theorem~\ref{koenigs-thm} is proved by defining a power series and
showing that it converges.  This approach is superficially
constructive, and can be used to obtain numerical approximations to
values of $\phi$.  It is also essentially the same as
the construction we applied in Section~\ref{find phi sxn}.  But unless
the constructed $\phi$ can be expressed in terms of familiar
functions, this approach is also unlikely to produce exact expressions for the limits
of candidate sequences.  

As the preceding remarks show, the ideas of candidate sequences and
Koenigs functions are closely intertwined.  This connection does not
seem to help much in finding exact limits of candidate sequences, but
it can be used in another way to show that candidate sequences
converge properly.  To illustrate this point, we provide the following
example.

\paragraph{Candidate Sequences for a Family of Root-Like Functions.}

We consider the family of functions 

\[ {\cal F} = 
\left\{ f_L(t) = \sqrt{L^2 - L + t} \;\;\;\rule[-2.1ex]{.6pt}{6ex}
\;\; L > 1. \right\}  \]

Each $f_L$ in the family is root-like, with an attracting fixed point at
$L$, and with $m_L = f_L'(L) = 1/(2L) \in (0,1).$  These functions are
discussed in~\cite{curriefunc}, where they are used to define a
function $C$.  As a natural extension of the mysterious pattern,
$C(L)$ is defined in terms of a modified candidate sequence for $f_L$
with an initial value $t_0 = 0$.  Specifically
\begin{equation} \label{C(L)def}
 C(L) = \sqrt{\lim_{n\rightarrow \infty} (2L)^n (L - f_L^{(n-1)}(0))}. 
\end{equation}
In~\cite{curriefunc} it is shown by elementary means that the limit
defining C(L) converges for $L>1$.

Convergence also follows from our Theorem~\ref{thm1}, because every
$f_L$ is root-like in an interval $I$ containing $L$.  If $0$ is in
that interval then Theorem~\ref{thm1} shows directly that $C(L)$ is
defined.  But even if that is not the case, the sequence $t_n =
f_L^{(n)}(0)$ increases torward $L$.  Thus, for some $k \ge 0$
$f^{(k)}(0) \in I$, and the candidate sequence with 
$t_0 = f^{(k)}(0)$ converges
properly by Theorem~\ref{thm1}.  Then the candidate
sequence with $t_0=0$ also converges properly by the same
logic used to prove Lemma~\ref{lemma3}.

We proceed now to present a third proof of convergence, as an
application of Theorem~\ref{koenigs-thm}.  The argument involves a
single value of $L$, so for simplicity we express $f_L$ simply as $f$.

As a first step we extend $f$ into the complex domain.  Specifically,
we consider $f(z) = \sqrt{L^2-L+z}$ for $z$ in 
$D_L(L)$, the open disk of radius $L$ centered at $L$ on the real
axis in $\bbC$. 
Then $\{L^2 - L + z \;|\; z \in D_L(L)\}$ is simply the translation of
$D_L(L)$ by $L^2-L$, and so is the
congruent disk centered at $L^2$, namely $D_L(L^2)$.  This is disjoint
from the negative real axis in $\bbC$ 
so the (principle value of the) complex square root is a
holomorphic function there.  In polar form, if $z = re^{i\theta}$ with
$r \ge 0$ and $-\pi/2 \le \theta \le \pi/2$,
then $\sqrt{z} = \sqrt{r}\cdot e^{i\theta / 2}$.  Thus we can envision
the square root as the polar coordinates map $(r,\theta) \mapsto
(\sqrt{r},\theta/2)$ in the real plane.  This shows that the
squareroot map is injective, hence $f$ is as well.  It is also the
case that $f$ maps $D_L(L)$ into itself.

Now the action of $f$ on $D_L(L)$ can be shifted to the open unit disk $D$
by means of the map $\psi(z) = L+Lz$ and its inverse 
$\psi^{-1}(z) = (z-L)/L$.  We define
\[ g(z) = \psi^{-1}(f(\psi(z))).   \]
This is a holomorphic map of $D$ into itself, with a fixed point at the
center $0$.  The multiplier $\lambda = g'(0)$ can be computed from the
equation 
\begin{equation} \label{psiog=fopsi}
 \psi(g(z)) = f(\psi(z)).
\end{equation}
Differentiating,
$\psi'(g(z))g'(z) = f'(\psi(z))\psi'(z)$.  Since $\psi'(z)=L$, at $z = 0$ we obtain
$Lg'(0) = f'(L)L$.  Thus $\lambda = f'(L) = 1/(2L) \in (0,1)$ showing
that 0 is an attracting fixed point of $g$.  As usual, let $m = f'(L)
= \lambda$.

Let $\phi$ be the holomorphic function defined on a disk $D_\epsilon(0)$ whose existence is
asserted by Theorem~\ref{koenigs-thm}, such that
$\phi(g(z))=\lambda\phi(z).$   This implies that
$\phi(g^{(n)}(z)) = \lambda^n \phi(z))$.  And because $\phi'(0) = 1,$
we can choose $\epsilon$ such that $\phi$ is injective on $D_\epsilon(0)$.

We know that for $z \in D_\epsilon(0)$,
$\phi(z) = \lim_{n \rightarrow \infty} g^{(n)}(z)/m^n$.  
And $w = \psi(z) \in D_{\epsilon L}(L)$. But
\[ g^{(n)}(z) =  \psi^{-1} f^{(n)}\psi(z) 
= \psi^{-1} f^{(n)}(w) \]
so 
\[ \psi g^{(n)}(z) =  Lg^{(n)}(z)+L = f^{(n)}(w).\] 
Thus $Lg^{(n)}(z) = f^{(n)}(w) - L $  and
\[  L\frac{g^{(n)}(z)}{m^n} = \frac{f^{(n)}(w) - L}{m^n}.   \]
The convergence of the sequence on the left, implies the
convergence of the sequence on the right.  The limit on the left
vanishes at $z=0$, giving $\phi(0) = 0.$ But by injectivity, $\phi(z)
\ne 0$ for $z \ne 0.$  Therefore the limit on the right converges properly
except when $z = 0$ and thus when $w = L.$  In conclusion, we have
shown that when $w \in D_{\epsilon L}(L)$ the corresponding candidate
sequence converges, and for $w \ne L$ it converges properly.
Furthermore, as noted earlier, this same conclusion must hold for $w$
in the basin of attraction of $f,$ which includes $D_L(L)$.

Recall that our ultimate goal is to establish convergence in the definition of
$C(L)$. To that end, we rewrite~\eqref{C(L)def} in the form
\[ C(L) = \sqrt{4L^2\lim_{n\rightarrow \infty} 
           (2L)^{n-2} (L - f^{(n-2)}(f(0)))}. \]
This version contains the limit of a candidate sequence $c_n$ for $f$ with
$t_0 = f(0) = \sqrt{L^2 - L}.$  We know $L>1,$ so $0 < t_0 < L$ which
shows that $t_0 \in D_L(L)$.  As we have seen, the candidate sequence
for this $t_0$ therefore converges to a nonzero limit, and $C(L)$ does
exist.

The point of this example has been to show that the Koenigs function
approach can be applied to demonstrate proper convergence of candidate
sequences.  In the process, we have seen how a change of variables
can shift the context of Theorem~\ref{koenigs-thm} to the case
of an injective holomorphic function $f$ on an arbitrary disk $G$
centered at a fixed point $L$.  For a real function $f$ as discussed
in Section~\ref{eigfunc sxn}, this may provide a way to establish proper convergence of
candidate sequences.  In order for this to be possible, the real
function $f$ must extend to a holomorphic function in a complex neighborhood of an attracting
fixed point.

\begin{figure}[hbt]
\centering 
\includegraphics[scale=0.5]{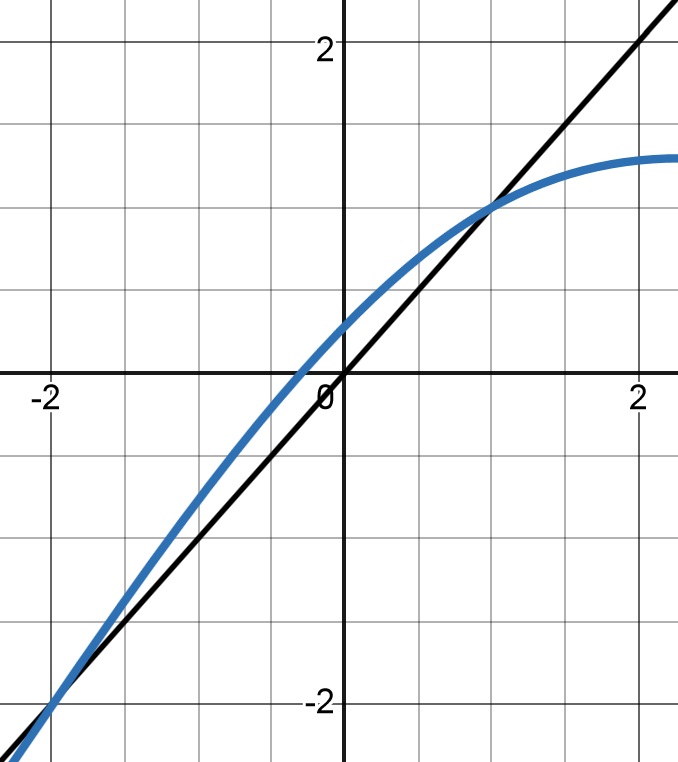}
\caption{Graph of $f$ near $(1,1)$.  Consistent with this graph it can
be shown that $f$ is root-like between $(-7+\sqrt{37.5})/3\approx -.29$ and 
$(13-\sqrt{37.5})/3 \approx 2.3$.}

\label{desmos-graph-4} 
\end{figure}

\paragraph{Comparing Theorem~\ref{thm1} and Theorem~\ref{koenigs-thm}.}
How does the use of Theorem~\ref{koenigs-thm} to demonstrate proper
convergence of candidate sequences compare with direct application of
Theorem~\ref{thm1}? 
It is not difficult to construct a function $f$ for which
Theorem~\ref{thm1} can be applied, but not Theorem~\ref{koenigs-thm}.  
For example, with $g(t) = .5t -.24t^2+.024t^2 |t|$ we see that 
$g(0) = 0$, $g'(0) = .5$, and $g''(0) = -.48$.  Thus $g$ is increasing
and concave down near the fixed point at $t=0.$  
Now define $f(t) = 1 + g(t-1)$ taking the fixed point of $g$ to
the point $(1,1)$.  See~\figref{desmos-graph-4}.
The graph of $f$, like that of $g$, is increasing and concave down
near the fixed point, and so is root-like there.  As the graph shows,
$f$ is root-like in the interval $(-.28,2.2)$.  For $t_0$ in that same
interval, Theorem~\ref{thm1} shows that candidate sequences for $f$ converge properly.

However, we cannot apply Theorem~\ref{koenigs-thm} to demonstrate this
result.  This follows because  
$g''(t) = .144|t|-.48$ hence $g'''(0)$ is undefined.  Therefore
$f'''(1)$ is also undefined, and so $f$ cannot be extended
to a holomorphic function in a neighborhood of the fixed point.

On the other hand there exist real functions that are not root-like,
but for which extension to the complex setting does imply properly
convergent candidate sequences.  For example, take $g(z) = z/2 -
z^4/12,$ which satisfies the hypotheses of Theorem~\ref{koenigs-thm},
and define $f(t) = 1 + g(t-1)$.  Then $f''(1) =g''(0)=0$ so $f$ is not
root-like at 1, its fixed point, and so Theorem~\ref{thm1} does not
apply.  But a closer examination shows that $f'$ is strictly
decreasing everywhere, and hence concave down, inspite of the fact
that $f''(1)=0$.  This suggests that the definition of root-like
underlying Theorem~\ref{thm1} can be relaxed.  If so, proper
convergence of candidate sequences for this example may be implied by
both Theorem~\ref{thm1} and Theorem~\ref{koenigs-thm}.

In terms of the formulation of Theorem~\ref{thm1} in this paper, 
we see that neither theorem supercedes the other as regards
demonstrating proper convergence of candidate sequences.  That is,
there are examples where Theorem~\ref{thm1} applies but not
Theorem~\ref{koenigs-thm}, and vice versa.  It may be
that Theorem~\ref{thm1} can be extended to cover all cases where
Theorem~\ref{koenigs-thm} applies, but that remains to be seen.

In this context, we might also observe that we have witnessed a successful
application of the eigen-function approach that could not arise from
Theorem~\ref{koenigs-thm}.  In particular we encountered the eigen-function
$\phi(\theta) = 2\cos(\theta)$ for $f(t) = \sqrt{2+t}$ in
Section~\ref{eigfunc sxn}.  For that example, $\phi'(0) = 0$, which is
incompatible with the assumptions of Theorem~\ref{koenigs-thm}.  It is
possible that an analogue of Theorem~\ref{koenigs-thm} that applies in
the case of super-attracting fixed points can be used here, but we
have not yet explored this idea.  

In this same example we also found an alternate eigen-function
$\phi(\sqrt{-\theta})$ whose derivative does not vanish at 0.  It
remains to be seen whether this foreshadows a more general development
concerning contractions $f$ admitting multiple non-equivalent $\phi$
functions, not all of which exhibit the condition $\phi'(0) = 0$. 

\paragraph{Verifying A Final Condition on $f$.}  Earlier we
indicated that
$f(z) = \sqrt{L^2-L+z}$ maps $D_L(L)$ into itself.  For completeness
sake, we now prove this assertion.

\begin{figure}[hbt]
\centering 
\includegraphics[scale=0.5]{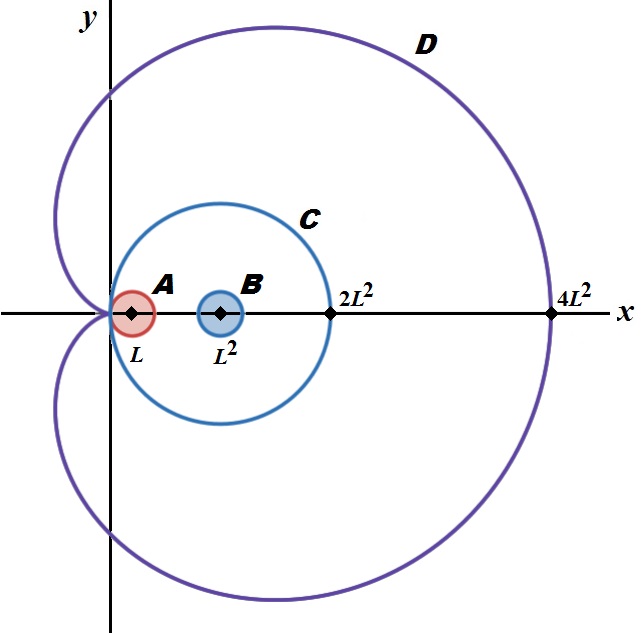}
\caption{Translating $A = D_L(L)$ by $L^2-L$ produces 
$B=D_{L}(L^2) \subset D_{L^2}(L^2) = C.$ In turn, $C$ is contained in
the cardioid $D,$ which is $S(A)$ where $S$ is the map $z \mapsto z^2.$}

\label{disk+cardoid} 
\end{figure}

We know $D_L(L)$ is a filled circle centered on the real axis at $L >
1$ and tangent to the imaginary axis at the origin.  This is shown
in~\figref{disk+cardoid} as the shaded circle $A$.  We will consider
the figure to represent both \bbC\ and \bbRR, identifying the
complex number $z = x+iy$ with the point $(x,y)$.

For $z \in D_L(L)$ let $w = L^2-L + z$.  This is just a translation of
$z$ by $L^2-L$ units to the right, and the image of $D_L(L)$ is thus
$D_L(L^2)$, shown in the figure as the shaded disk $B$.  Our goal is
to show that the image of $B$ under the square root map $Q$ is contained
in $A$.  However, this is equivalent to showing that $B$ is
contained in the image of $A$ under the square map $S$, because both
$Q$ and $S$ are injective on the closed halfplane $x \ge 0$.  And as
$B$ is contained in $D_{L^2}(L^2)$ (shown in the figure as the
unshaded disk $C$), it suffices to show that $C$ is
contained in $S(A)$.  To accomplish this, we show that $S(A)$ is the
region bounded by the cardioid $D$ in the figure.

Recall that the polar coordinates formula for the circle of radius $R$
centered at $R$ is $r = 2R\cos(\theta),$ for $\theta \in [-\pi/2,\pi/2]$.
Therefore using polar coordinates  we can describe $D_L(L)$ as 
\[A =\{ (r,\theta)\;|\; -\pi/2 \le \theta \le \pi/2, 0 \le r \le 2L\cos(\theta)\}\]
and $D_{L^2}(L^2)$ as
\[C =\{ (r,\theta)\;|\; -\pi/2 \le \theta \le \pi/2, 0 \le r \le 2L^2\cos(\theta).\}\]

At a particular point $z_0 = (r_0, \theta_0) \in A$ we have
$S(z_0) = (r_0^2,2\theta_0)$.  The polar coordinates of this point
are $r_1 = r_0^2$ and $\theta_1 = 2\theta_0$.  Now observe that
\[ 0 \le r_0 \le 2L\cos(\theta_0) = 2L\cos(\theta_1/2).    \] 
Squaring both sides produces
\[ 0 \le r_0^2 = r_1  \le  4L^2\cos^2(\theta_1/2)=2L^2(1+\cos(\theta_1)).  \] 
Thus, the polar equation of the boundary curve for $S(A)$ is
\[ r = 2L^2(1+\cos(\theta))  \]
which we recognize as the cardioid shown in the figure.

Visually, it is apparent that $C$ is contained in $D$.  Analytically,
at an arbitrary point $z = (r,\theta) \in C$ we have 
$0 \le r \le 2L^2\cos(\theta) \le 2L^2(1+\cos(\theta)) $.  Thus $z \in
D,$ which proves that $C$ is contained in $S(A)$ as desired. 

It is now easy to see that $f$ does not map $D_L(L)$ {\em onto}
itself.  If it did, the squareroot function would have 
to map $B$ onto $A,$ and that would imply that $S$ maps $A$ into $B$.
This is false because $z = L+2/3 < 2L$ but $z^2 = L^2 + 4/3L + 4/9 >
L^2+L$.  The first inequality says that $z \in A$ and the second says
that $S(z) \not\in B.$

This completes our consideration of Koenigs functions.  Our final
topic involves some additional analogs of the mysterious pattern.
These are obtained using trigonometric identities for $\cos(\theta/K)$
guided by the use of half-angle formulas in the original mysterious
pattern.  Chebyshev polynomials arise naturally in this development.

\section{Inverse Chebyshev Polynomials} \label{Cheby_stuff}

The goal of the development to follow is a construction that
produces limit values of the form
\[ \lim_{n \rightarrow \infty} c_n = \pi^2/8   \]
where
\[ c_n = (K^2)^n  \left| 1 - f^{(n)}(0) \right|    \]
and where $n$ is the number of nested applications of $f$ and $K$ is a
positive integer.  Here, $c_n$ is roughly equivalent to the square of
the mysterious pattern's $a_n$.

As a preliminary step, let us see that the case with $K=2$ is
equivalent to the limit in mysterious pattern.  For this case
\[ f(t) = \sqrt{\frac12 + \frac{t}{2}} = \frac12\sqrt{2+2t}. \]
When we apply $f$ several times to a starting value of 0, we generate
the following:
\begin{eqnarray*}
f(0)   & = & \frac12\sqrt{2}   \\
f^{(2)}(0)   & = & \frac12 \sqrt{2+\sqrt{2}}   \\
f^{(3)}(0)   & = & \frac12 \sqrt{2+\sqrt{2+\sqrt{2}}}.
\end{eqnarray*}
These are almost identical to the nested radical sequence that appears in the
original mysterious pattern, except that
each term has an extra factor of $1/2$.  Therefore the 
limit of $f^{(n)}(0)$ will be $1$ instead of $2$.  So,
\[ c_n = 4^n \left| 1 - \frac12 \sqrt{2+\sqrt{2+\cdots+\sqrt{2}}} \;\right| \]
and
\[ 2c_n = 4^n \left| 2 - \sqrt{2+\sqrt{2+\cdots+\sqrt{2}}} \;\right|. \]
Here we see that $\sqrt{2c_n}$ differs from $a_n$ only in the fact that
the former has $n+1$ nested radicals instead of the latter's $n$
radicals.  To correct this, multiply by 4 to obtain
\[ 8c_n = 4^{n+1} \left| 2 - \sqrt{2+\sqrt{2+\cdots+\sqrt{2}}} \;\right| \]
and now we see that $\sqrt{8c_n} = a_{n+1}$.  The starting assertion
of this discussion, that the limit of $c_n$ is
$\pi^2/8,$ we now observe is equivalent to $8c_n$ having a limit of
$\pi^2$.  Thus, as asserted, the case with $K = 2$ is equivalent to
the mysterious pattern.

As we have seen the original mysterious pattern can be understood in
terms of the eigen-function $\phi(t) = 2\cos(t)$.  In particular we
noticed that the iterated function in that case, $f(t) = \sqrt{2+t}$
could be expressed in the form $\phi \comp \mu_{1/2} \comp
\phi^{-1}(t)$, or equivalently $f(\phi(\theta)) = \phi(\frac12
\theta)$.  This follows from the half-angle formula for cosine, which
in turn can be derived from the double angle formula.  Mimicing this
structure, we wish to formulate an identity for $\cos(\theta/K)$
starting from an identity for $\cos(K\theta)$.

\subsection{Cosine of $K$ Theta}  It is well known that for integer
$K$, $\cos(K\theta)$ can be expressed as a polynomial in $\cos
\theta$.  One way to see this is to use the identity 
$e^{iK\theta} = \left(e^{i\theta}\right)^K$.  Thus
\[\cos(K\theta) + i\sin(K\theta) 
    = (\cos(\theta) + i\sin(\theta))^K. \]

The right side can be expanded using the binomial theorem, and then
separated into real and imaginary parts.  Any factors of
$\sin(\theta)$ that show up in the real part have to have even exponents,
and so can be converted into expressions in $\cos^2(\theta)$.  Thus
$\cos(K \theta)$ can be expressed in the form $p_K(\cos \theta)$ for a
polynomial $p_K.$  These are in fact Chebyshev polynomials, but that
is not really relevant in what will follow.  The first few $p_K$ are
shown below.  Their graphs are displayed in~\figref{p_k graphs}.

\[ \begin{array}{|c|l|} \hline
 K & p_K(t) \\ \hline
 2 & 2t^2 - 1 \\
 3 & 4t^3 - 3t \\
 4 & 8 t^4 - 8 t^2 + 1\\
 5 & 16t^{5}-20t^{3}+5t \\ \hline
   \end{array}  \]

\begin{figure}[hbt]
\centering 
\includegraphics[width=4.5in]{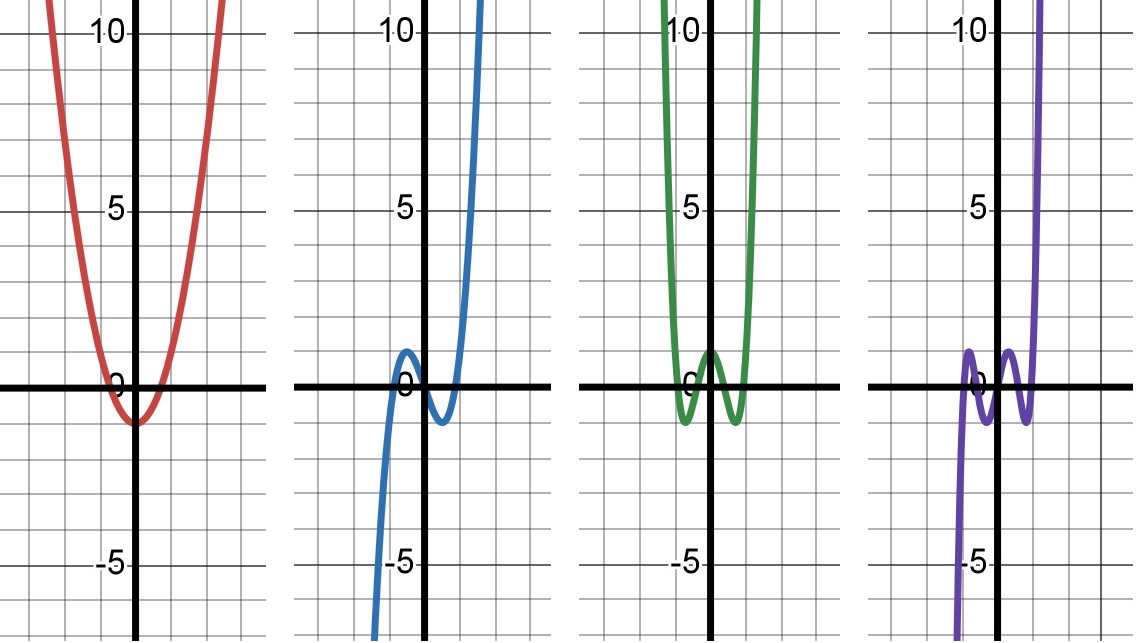}
\caption{Graphs of $p_K$ for $K = 2, 3, 4, 5,$ in order from left to right.
Notice that local minima on each graph occur with $p_K(t) = -1.$  
From the right-most local minimum on each graph the curve increases.}
\label{p_k graphs}
\end{figure}

To relate this to $\cos(\theta/K)$ we need the inverse function
of $p_K$.  For $K = 2$ the inverse of $p_K$ is 
\[ f_2(t) = \sqrt{\frac{1+t}{2}} = \frac12 \sqrt{2+2t} \]
and we have already seen how this function is connected with the mysterious pattern.

\begin{figure}[hbt]
\centering 
\includegraphics[width=3in]{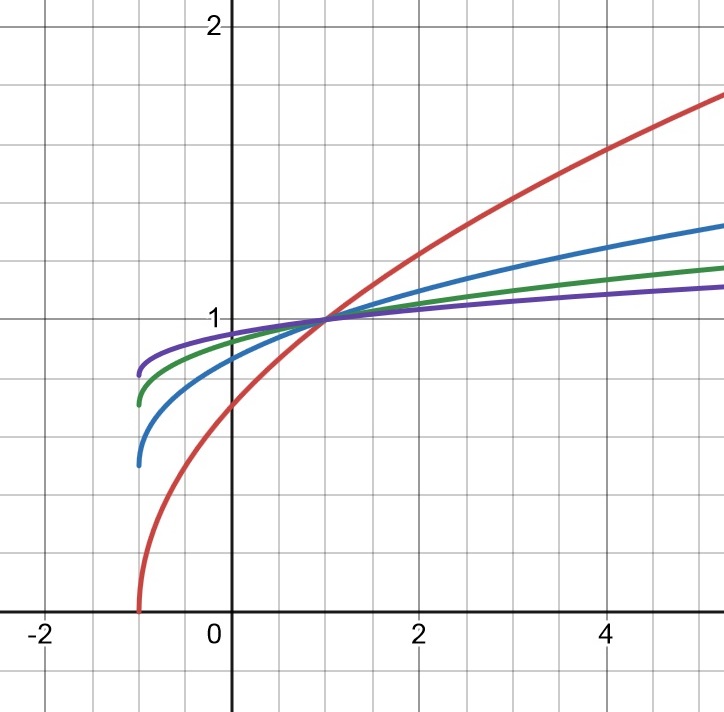}
\caption{Graphs of $f_K$ for $K = 2, 3, 4, 5.$  At the left end points
of the curves, they appear in order from bottom to top.  That is, at the
lowest left endpoint $K = 2,$ at the second lowest $K = 3,$ and so
on.  All of these curves have the same general shape as $f_2:$
increasing without bound and concave down, with a vertical tangent at
the left end-point.  For each curve $(1,1)$ represents a fixed
point.}
\label{f_k graphs} 
\end{figure}

More generally, we will define $f_K$ as the inverse function of $p_K$
on the maximal interval of the form $[z,\infty)$ for which $p_K$ is
increasing.  Visually, this corresponds to the part of the graph to the
right of the right-most local minimum.  As these minima all occur with
$p_K(z) = -1$, the inverse functions each have domain $[-1,\infty)$.
The graphs of $f_K$ for $K = 2, 3, 4, 5$ appear 
in~\figref{f_k graphs}.

From the way $p_K$ is 
defined, we know the identity 
\begin{equation} \label{mel.eq.1} 
\cos(K\theta) = p_K(\cos(\theta)).
\end{equation}
Taking $\theta=0$ shows that $p_K(1) = 1.$  We will invert $p_K$ in a
neighborhood of this fixed point.  It can be shown, using properties
of Chebyshev polynomials, that $p_K'(1) = K^2$ for $K \ge 1$.  Thus
$p_K$ is increasing in a neighborhood of 1, and so invertible on that
neighborhood.  The inverse of $p_K$ on any such neighborhood of 1 is
what we will denote as $f_K$.

For $-1 \le t \le 1$ the function $f_K$ is given by
\begin{equation} \label{f cos eqn}
  f_K(t) = \cos \left(\frac{\cos^{-1}(t)}{K} \right).
\end{equation}
This equation can be derived from~\eqref{mel.eq.1} by applying $f_K$ on both
sides and defining $t = \cos(K\theta)$.  It shows that \; 
$f_K : [-1,1] \rightarrow [\cos(\pi/K),1]$.\; 
For $t \ge 1$ $f_K$ can be represented by the equation
\begin{equation} \label{f cosh eqn}
  f_K(t) = \cosh \left(\frac{\cosh^{-1}(t)}{K} \right). 
\end{equation}
For our present purposes, we will consider only the case for $t < 1.$

The two preceding equations also show that
\[ f_K = \phi \comp \mu_{1/K} \comp \phi^{-1}  \]
where $\phi(\theta)$ equals $\cos(\theta)$ in the context of~\eqref{f cos eqn}
and $\cosh(\theta)$ for~\eqref{f cosh eqn}.  In either case $\mu_{1/K}(\theta) =
\theta/K.$  Thus we recognize that $\phi$ is an eigen-function for $f_K$.

At the same time, we know that $f_K$ is smooth because it is the
inverse function of a polynomial.  We also know that 
$f_K$ has a fixed point at $1$ and that $f_K'(1) = 1/K^2 < 1.$  This
implies that $f_K$ is a contraction near $1$.  In
fact, as suggested by the graphs of $f_K$ and $p_K$, $f_K'$ decreases
for $t \in (-1,\infty)$ and must equal 1 at a unique point $z < 1$ in that
interval.  Then $f_K'(t) \in (0,1)$ for all $t>z$ showing that $f$ is
a contraction on $(z,\infty)$.  For any $t_0$ in this interval
$f^{(n)}(t_0)$ converges monotonically to the fixed point $1.$  Thus it
makes sense to consider candidate sequences of the form
\[ c_n = \frac{|1 - f_K^{(n)}(t_0)|}{1/K^{2n}}.   \]
With the restriction $t_0 < 1$ we know that $f_K$ is defined by~\eqref{f cos eqn}.

\paragraph{Limits of Candidate Sequences.}
We have seen earlier that when $f = \phi \comp \mu_{\alpha} \comp
\phi^{-1}$ is a contraction the fixed point is $L = \phi(0)$.  In the
current context, where $\phi(\theta) = \cos(\theta),$ this reconfirms that there is a fixed point of $f_K$
at $t=1.$ However, we also observe that $\phi'(0) = 0$ and $\phi''(0)
= -1 \ne 0$.  Thus we can apply~\eqref{phi'(0)=0 eqn1} to find
\[ \lim_{n\rightarrow \infty} c_n 
   = \left| \frac{\theta_0^2 \phi''(0)}{2} \right|
   = \frac{\theta_0^2}{2}. \]
where $\cos(\theta_0) = t_0$.
As a particular instance, let $\theta_0 = \pi/2$ so that $t_0 = 0.$
Then
\[ \lim_{n \rightarrow \infty} c_n = \pi^2/8  .\]

\subsection{A New Radical Notation}
For each value of $K > 2$ we now have a new analog of the mysterious
pattern, complete with an exact expression for the limit.  Although
each $K$ produces the same limit, these analogs do not appear to be
equivalent to the original mysterious pattern, or to each other,
because the functions $f$ are distinct, as are the factors $(2K)^n$.

To make the analogy between the mysterious pattern and these new
examples more visual, let us introduce a new
algebraic symbol for $f_K$, writing $\sqrt[\langle K \rangle]{t} =
f_K(t-1)$.  Then $\sqrt[\langle K \rangle]{1+t} = f_K(t)$.  This is geometrically
reasonable because $f_K$ has a graph that resembles a $K$th-root curve
with its left endpoint at $t=-1$.  Algebraically, the function
$\sqrt[K]{t}$ is the inverse function of a particular $K$th degree
polynomial, whereas $\sqrt[\langle K \rangle]{t}$ is the inverse
function of a different particular $K$th degree polynomial.  Is it too
great a stretch to propose a notation for the latter that is similar
to the notation for the former?

We say it is not.  Accordingly, we now write the limit
of the candidate sequence as

\begin{equation} \label{ChebyLim1}
 \lim_{n \rightarrow \infty} 
   K^{2n}\left(1-\sqrt[\langle K \rangle]{1+
            \sqrt[\langle K \rangle]{1+ \cdots
            \sqrt[\langle K \rangle]{1}}}\right) = \pi^2/8. 
\end{equation}

As a specific example, take $K = 5.$  Then
\[ \lim_{n \rightarrow \infty} 
   25^{n}\left(1-\sqrt[\langle 5 \rangle]{1+
            \sqrt[\langle 5 \rangle]{1+ \cdots
            \sqrt[\langle 5 \rangle]{1}}}\right) = \pi^2/8. \]
Expressed in this way, these new examples stand as obvious analogs of
the mysterious pattern.

As one additional note, we can define a second eigen-function,
\[ \phi_2(\theta) = \cos(\sqrt{-\theta})  \]
for $\theta < 0.$  Then, echoing what we saw for the original
mysterious pattern at the end of Section~\ref{origMP sxn}, we have
$\phi_2 = \phi \comp \eta$ where $\eta(\theta) = 
\sqrt{-\theta}$.  Applying identical reasoning in the current context,
we find that 
\[\phi_2 \comp \mu_{1/k^2} \comp \phi_2^{-1} =
\phi \comp \mu_{1/k} \comp \phi^{-1} = f_K. \]
This confirms that $\phi_2$ is indeed an eigen-function and can be
used to compute limits of candidate sequences for $f_K.$  Also as
before, $\phi_2'(0) \ne 0.$  Working directly from the
definition of $\phi_2$, we can determine that $\phi_2'(0) = 1/2.$  Alternatively,
we can derive the power series for $\phi_2$ by substituting
$\sqrt{-\theta}$ for $t$ in the series for $\cos(t)$, thus obtaining
\[ \phi_2(\theta) = \sum_{k=0}^{\infty} \frac{\theta^k}{(2k)!}.  \]
The linear coefficient then gives $\phi_2'(0) = 1/2$.

In any case, knowing that $\phi_2'(0) \ne 0$ allows us to use 
use~\eqref{c_n lim eq2} to find the limit of a candidate sequence
$c_n$.  Thus,  for any valid $t_0,$ and with $\phi_2(\theta_0) = t_0,$
we  conclude
\[ \lim_{n \rightarrow \infty} c_n = |\phi_2'(0)\theta_0| =
|\theta_0|/2. \]
In particular, with $\theta_0 = -\pi^2/4,$ corresponding to $t_0 = 0$,
the preceding equation becomes
\[ \lim_{n \rightarrow \infty} c_n = \pi^2/8, \]
reconfirming the result we found using $\phi.$

Finally, it should be noted that the derivation of~\eqref{ChebyLim1}
involved some arbitrary choices, and somewhat different looking
results can be derived given different choices.  For example, in~\cite{geom_ext_ppr} 
the following variant for the $K=3$ case is derived:  
\[ \lim_{n \rightarrow \infty} 
3^{n} \sqrt{3-\sqrt[\blacktriangle]
      {3 + \sqrt[\blacktriangle]
      {3 + \cdots + \sqrt[\blacktriangle]{3}}}} 
= \frac{3\pi\sqrt{6}}{4}. \]
In that construction
\[ f(t) = 
3 \cos\left( \frac{1}{3}\cos^{-1}\left(\frac{t}{3}\right)\right)   \]
and
\[ \sqrt[\blacktriangle]{t} =  f(t-3).  \]

The analysis leading to~\eqref{ChebyLim1} can be modified in a similar
fashion by defining
$\phi(t) = K\cos(t)$, so that the expression 
$\phi \mu_k \phi^{-1}$ leads to $f(t) = K p_K^{-1}(t/K)$ and the fixed
point occurs at $t=K$.  We can follow the identical development as
before, ultimately defining 
$\sqrt[\langle K \rangle]{t} = f_K(t-K)$.  Then in place
of~\eqref{ChebyLim1} we arrive at 
\begin{equation} \label{ChebyLim2}
 \lim_{n \rightarrow \infty} 
   K^{2n}\left(K-\sqrt[\langle K \rangle]{K+
            \sqrt[\langle K \rangle]{K+ \cdots
            \sqrt[\langle K \rangle]{K}}}\right) = K\pi^2/8. 
\end{equation}
We leave it to the reader to determine what variation on these results
best deserves to be considered a novel analog of the original
myaterious pattern.

\section*{Conclusion}
We have seen in this paper that Currie's mysterious
pattern can be generalized to various instances of properly converging
candidate sequences for a contraction in a neighborhood of its unique
fixed point.  For the set of functions we have described as root-like,
we found that candidate sequences always converge properly, though we
did not in general determine the limits.  The convergence result was
proved by approximating the root-like function $f$ with a root-like
M\"{obius} function $Q$.  The convergence of candidate sequences for
such a $Q$ implies the convergence of candidate sequences for the
corresponding $f$.  Moreover, we discovered exact values of limits of
candidate sequences for all such $Q$, thus demonstrating proper
convergence.

We also analyzed candidate sequences in the case of a function $f$
expressible as a composition $\phi \mu_{\alpha} \phi^{-1}$.  This
gives rise to results reminiscent of matrix diagonalization, and in
particular leads to candidate sequence limits expressed in terms of
$\phi.$  Various ways to identify $f$-$\phi$ pairs were explored, though
without leading to new exact limits of candidate sequences.  However,
this development led to a connection with complex analysis and Koenigs
functions.  It also enabled us to find a collection of new analogs for
the mysterious pattern, including exact values of limits of candidate
sequences, based on identities for $\cos(\theta/K)$.  For these
examples we defined $f$ by inverting the $K$th Chebyshev polynomial on
part of its domain, and identified $\cos(\theta)$ as the corresponding
$\phi$ function.  By introducing a new radical notation to represent
these $f$'s, a nice visual analog of the mysterious pattern was obtained.

\end{document}